\documentclass{article}
\usepackage{amssymb}

\input{tcilatex}
\begin{document}

\begin{center}
\bigskip \textbf{Hardy's uncertainty principle} \textbf{and unique
continuation properties for abstract Schr\"{o}dinger equations}

\textbf{Veli Shakhmurov}

Department of Mechanical Engineering, Okan University, Akfirat, Tuzla 34959
Istanbul, Turkey,

E-mail: veli.sahmurov@okan.edu.tr

A\textbf{bstract}
\end{center}

In this paper, Hardy's uncertainty principle and unique continuation
properties of abstract Schr\"{o}dinger equations in vector-valued $L^{2}$
classes are obtained.


\begin{center}
\bigskip\ \ \textbf{AMS 2010: 35Q41, 35K15, 47B25, 47Dxx, 46E40 }

\textbf{1. Introduction, definitions}
\end{center}

\bigskip In this paper, the unique continuation properties of the abstract
Schr\"{o}dinger equations

\begin{equation}
i\partial _{t}u+\Delta u+Au+V\left( x,t\right) u=0,\text{ }x\in R^{n},\text{ 
}t\in \left[ 0,T\right] ,  \tag{1.1}
\end{equation}%
are studied, where $A$ is a linear operator$,$ $V\left( x,t\right) $ is a
given potential operator function in a Hilbert space $H$, subscript $t$
indicates the partial derivative with respect to $t$, $n$ is the dimension
of space variable $x$, $\Delta $ denotes the Laplace operator in $R^{n}$ and 
$u=$ $u(x,t)$ is the $H$-valued unknown function. This linear result was
then applied to show that two regular solutions $u_{1}$ and $u_{2}$ of
non-linear Schr\"{o}dinger equations 
\begin{equation}
i\partial _{t}u+\Delta u+Au=F\left( u,\bar{u}\right) ,\text{ }x\in R^{n},%
\text{ }t\in \left[ 0,T\right] ,  \tag{1.2}
\end{equation}%
and for very general non-linearities $F$, must agree in $R^{n}\times \lbrack
0,T]$, when $u_{1}-u_{2}$ and its gradient decay faster than any quadratic
exponential at times $0$ and $T$.

Hardy's uncertainty principle and unique continuation properties for Schr%
\"{o}dinger equations studied e.g in $\left[ \text{4-7}\right] $ and the
referances therein.\ In contrast to the mentioned above results we will
study the unique continuation properties of abstract Schr\"{o}dinger
equations with operator potentials. Abstract differential equations studied
e.g. in $\left[ \text{2, 12-19, 22, 24, 25}\right] .$ Since the Hilbert
space $H$ is arbitrary and $A$ is a possible linear operator, by choosing $H$
and $A$ we can obtain numerous classes of Schr\"{o}dinger type equations and
its systems which occur in a wide variety of physical systems. Our main goal
is to obtain sufficient conditions on a solution $u$, the operator $A,$
potential $V$ and the behavior of the solution at two different times, $%
t_{0} $ and $t_{1}$ which guarantee that $u\left( x,t\right) \equiv 0$ for $%
x\in R^{n},$ $t\in \left[ 0,T\right] $. If we\ choose the abstract space $H$
a concrete Hilbert space, for example $H=L^{2}\left( \Omega \right) $, $A=L,$
where $\Omega $ is a domin in $R^{m}$ with sufficientli smooth boundary and $%
L$ is elliptic operator then, we obtain the unique continuation properties
of followinng Schr\"{o}dinger equation%
\begin{equation}
\partial _{t}u=i\left[ \Delta u+Lu+V\left( x,t\right) u\right] ,\text{ }x\in
R^{n},\text{ }y\in \Omega ,\text{ }t\in \left[ 0,T\right] .  \tag{1.3}
\end{equation}

\ Moreover, let we choose $H=L^{2}\left( 0,1\right) $ and $A$ to be
differential operator with generalized Wentzell-Robin boundary condition
defined by 
\begin{equation}
D\left( A\right) =\left\{ u\in W^{2,2}\left( 0,1\right) ,\text{ }%
B_{j}u=Au\left( j\right) +\dsum\limits_{i=0}^{1}\alpha _{ij}u^{\left(
i\right) }\left( j\right) ,\text{ }j=0,1\right\} ,\text{ }  \tag{1.4}
\end{equation}%
\[
\text{ }Au=au^{\left( 2\right) }+bu^{\left( 1\right) }+cu, 
\]%
where $\alpha _{ij}$ are complex numbers, $a=a\left( y\right) ,$ $b=b\left(
y\right) $, $c=c\left( y\right) $ are complex-valued functions and $V\left(
x,t\right) $ is a integral operator defined by 
\[
V\left( x,y,t\right) u=\dint\limits_{0}^{1}K\left( x,y,\tau ,t\right)
u\left( x,y,\tau ,t\right) d\tau , 
\]%
where, $K=K\left( x,y,\tau ,t\right) $ is complex valued bounded function.
Then, we get the unique continuation properties of the Wentzell-Robin type
boundary value problem (BVP) for the following Schr\"{o}dinger equation 
\begin{equation}
\partial _{t}u=i\left[ \Delta u+a\frac{\partial ^{2}u}{\partial y^{2}}+b%
\frac{\partial u}{\partial y}+cu+\dint\limits_{0}^{1}K\left( x,y,\tau
,t\right) u\left( x,y,\tau ,t\right) d\tau \right] ,\text{ }  \tag{1.5}
\end{equation}%
\[
x\in R^{n},\text{ }y\in \left( 0,1\right) ,\text{ }t\in \left[ 0,T\right] , 
\]%
\ \ \ 

\begin{equation}
B_{j}u=Au\left( x,j,t\right) +\dsum\limits_{i=0}^{1}\alpha _{ij}u^{\left(
i\right) }\left( x,j,t\right) =0\text{, }j=0,1.  \tag{1.6}
\end{equation}

Note that, the regularity properties of Wentzell-Robin type BVP for elliptic
equations were studied e.g. in $\left[ \text{11, 12 }\right] $ and the
references therein. Moreover, if put $H=l_{2}$ and choose $A$ as a infinite
matrix $\left[ a_{mj}\right] $, $m,j=1,2,...,\infty ,$ then we obtain the
unique continuation properties of the following system of Schr\"{o}dinger
equation 
\begin{equation}
\partial _{t}u_{m}=i\left[ \Delta u_{m}+\sum\limits_{j=1}^{N}\left(
a_{mj}+b_{mj}\left( x,t\right) \right) u_{j}\right] ,\text{ }x\in R^{n},%
\text{ }t\in \left( 0,T\right) .  \tag{1.7}
\end{equation}

Let $E$ be a Banach space. $L^{p}\left( \Omega ;E\right) $ denotes the space
of strongly measurable $E$-valued functions that are defined on the
measurable subset $\Omega \subset R^{n}$ with the norm

\[
\left\Vert f\right\Vert _{L^{p}}=\left\Vert f\right\Vert _{L^{p}\left(
\Omega ;E\right) }=\left( \int\limits_{\Omega }\left\Vert f\left( x\right)
\right\Vert _{E}^{p}dx\right) ^{\frac{1}{p}},\text{ }1\leq p<\infty \ . 
\]

For $p=2$ and $H$ Hilbert space we get Hilbert space of $H$-valued functions
with inner product of two elements $f$, $g\in L^{2}\left( \Omega ;H\right) $:%
\[
\left( f,g\right) _{L^{2}\left( \Omega ;H\right) }=\int\limits_{\Omega
}\left( f\left( x\right) ,g\left( x\right) \right) _{H}dx. 
\]

Let $C\left( \Omega ;E\right) $ denote the space of $E-$valued, bounded
uniformly continious functions on $\Omega $ with norm 
\[
\left\Vert u\right\Vert _{C\left( \Omega ;E\right) }=\sup\limits_{x\in
\Omega }\left\Vert u\left( x\right) \right\Vert _{E}. 
\]

$C^{m}\left( \Omega ;E\right) $\ will denote the spaces of $E$-valued
bounded uniformly strongly continuous and $m$-times continuously
differentiable functions on $\Omega $ with norm 
\[
\left\Vert u\right\Vert _{C^{m}\left( \Omega ;E\right) }=\max\limits_{0\leq
\left\vert \alpha \right\vert \leq m}\sup\limits_{x\in \Omega }\left\Vert
D^{\alpha }u\left( x\right) \right\Vert _{E}. 
\]

$O_{R}=\left\{ x\in R^{n},\text{ }\left\vert x\right\vert <R\right\} $ for $%
R>0$. Let $\mathbb{N}$ denote the set of all natural numbers, $\mathbb{C}$
denote the set of all complex numbers.

Let $E_{1}$ and $E_{2}$ be two Banach spaces. $L\left( E_{1},E_{2}\right) $
will denote the space of all bounded linear operators from $E_{1}$ to $%
E_{2}. $ For $E_{1}=E_{2}=E$ it will be denoted by $L\left( E\right) .$

A linear operator\ $A$ is said to be positive in a Banach\ space $E$ with
bound $M>0$ if $D\left( A\right) $ is dense on $E$ and $\left\Vert \left(
A+sI\right) ^{-1}\right\Vert _{L\left( E\right) }\leq M\left\vert
s\right\vert ^{-1}$ for any $s\in $ $\left( -\infty ,0\right) ,$ where $I$
is the identity operator in $E.$

Let $\left[ A,B\right] $ be a commutator operator, i.e. 
\[
\left[ A,B\right] =AB-BA 
\]%
for linear operators $A$ and $B.$

Sometimes we use one and the same symbol $C$ without distinction in order to
denote positive constants which may differ from each other even in a single
context. When we want to specify the dependence of such a constant on a
parameter, say $\alpha $, we write $C_{\alpha }$.

\begin{center}
\textbf{2}. \textbf{Free} \textbf{absract Scr\"{o}dinger equation}
\end{center}

First of all, we generalize the result G. H. Hardy (see e.g $\left[ 20\right]
$, p.131) about uncertainty principle for Fourier transform:

\textbf{Lemma 2.1. }Let \textbf{\ }$f\left( x\right) $ be $H$-valued
function for $x\in R^{n}$ and%
\[
\left\Vert f\left( x\right) \right\Vert _{H}=O\left( e^{-\frac{\left\vert
x\right\vert ^{2}}{\beta ^{2}}}\right) ,\text{ }\left\Vert \hat{f}\left( \xi
\right) \right\Vert _{H}=O\left( e^{-\frac{4\left\vert \xi \right\vert ^{2}}{%
\alpha ^{2}}}\right) ,\text{ }x\text{, }\xi \in R^{n}\text{ for }\alpha
\beta <4. 
\]

Then $f\left( x\right) \equiv 0.$ Also, if $\alpha \beta =4$ then$\left\Vert
f\left( x\right) \right\Vert _{H}$ is a constant multiple of $e^{-\frac{%
\left\vert x\right\vert ^{2}}{\beta ^{2}}}.$

\textbf{Proof. }Indeed, by employing Phragmen--Lindel\"{o}f theorem to
Hilberts space valued analytic function class and by reasoning as in $\left[
8\right] $ we obtain the assertion.

Consider the Cauchy problem for free abstract Schr\"{o}dinger equation%
\begin{equation}
i\partial _{t}u+\Delta u+Au=0,\text{ }x\in R^{n},\text{ }t\in \left[ 0,1%
\right] ,  \tag{2.1}
\end{equation}

\begin{equation}
u\left( x,0\right) =f\left( x\right) .  \tag{2.2}
\end{equation}

The above result can be rewritten in terms of the solution of the $\left(
2.1\right) -\left( 2.2\right) $ on $R^{n}\times \left( 0,\infty \right) $ as:

Assume%
\[
\left\Vert u\left( x,0\right) \right\Vert _{H}=O\left( e^{-\frac{\left\vert
x\right\vert ^{2}}{\beta ^{2}}}\right) ,\text{ }\left\Vert u\left(
x,T\right) \right\Vert _{H}=O\left( e^{-\frac{\left\vert x\right\vert ^{2}}{%
\alpha ^{2}}}\right) \text{ for }\alpha \beta <4T. 
\]

Then $u\left( x,t\right) \equiv 0.$ Also, if $\alpha \beta =4T$, then $u$
has as a initial data a constant multiple of $e^{-\left( \frac{1}{\beta ^{2}}%
+\frac{i}{4T}\right) \left\vert y^{2}\right\vert }.$

\textbf{Lemma 2.2. }Assume $A$ is a pozitive operator in Hilbert space $H$
and $iA$ generates a semigrop $U\left( t\right) =e^{iAt}$. Then for $f\in
W^{s,2}\left( R^{n};H\right) $ there is a generalized solution of $\left(
2.1\right) $ expressing as%
\begin{equation}
u\left( x,t\right) =F^{-1}\left[ e^{iA_{\xi }t}\hat{f}\left( \xi \right) %
\right] ,\text{ }A_{\xi }=A+\left\vert \xi \right\vert ^{2}.  \tag{2.3}
\end{equation}

\textbf{Proof. }By applying the Fourier trasform to the problem $\left(
2.1\right) -\left( 2.2\right) $ we get 
\begin{equation}
i\partial _{t}\hat{u}\left( \xi ,t\right) +A_{\xi }\hat{u}\left( \xi
,t\right) =0,\text{ }x\in R^{n},\text{ }t\in \left[ 0,1\right] ,  \tag{2.4}
\end{equation}

\begin{equation}
\hat{u}\left( \xi ,0\right) =\hat{f}\left( \xi \right) ,\text{ }\xi \in
R^{n},  \tag{2.5}
\end{equation}

It is clear to see that the solution of the equation $\left( 2.4\right)
-\left( 2.5\right) $ can be exspressed as 
\[
\hat{u}\left( \xi ,t\right) =e^{iA_{\xi }t}\hat{f}\left( \xi \right) . 
\]

Hence, we obtain $\left( 2.3\right) .$

Let 
\[
X=L^{2}\left( R^{n};H\right) \text{, }Y^{k}=W^{k,2}\left( R^{n};H\right) 
\text{, }B=L^{\infty }\left( R^{n};L\left( H\right) \right) ,\text{ } 
\]%
\[
B=L^{\infty }\left( R^{n};B\left( H\right) \right) \text{ and }\mu \left(
t\right) =\frac{1}{\alpha \left( 1-t\right) +\beta t}. 
\]%
Consider the following abstract Schr\"{o}dinger equation%
\begin{equation}
\partial _{t}u=i\left[ \Delta u+Au+V\left( x,t\right) u\right] ,\text{ }x\in
R^{n},\text{ }t\in \left[ 0,1\right] ,  \tag{2.6}
\end{equation}
where $A$ is a linear operator in $H$ and $V\left( x,t\right) $ is a given
potential operator function in $H.$

Our main result in this paper is the following

\textbf{Theorem 1. }Assume the following condition are satisfied:

(1) $A$ is a symmetric operator in $H$ and $V\left( x,t\right) \in L\left(
H\right) $ for $\left( x,t\right) \in R^{n}\times \left[ 0,1\right] $;

(2) either, $V\left( x,t\right) =V_{1}\left( x\right) +V_{2}\left(
x,t\right) $, where $V_{1}\left( x\right) \in L\left( H\right) $ for $x\in
R^{n}$ and%
\[
\sup\limits_{t\in \left[ 0,1\right] }\left\Vert e^{\left\vert x\right\vert
^{2}\mu ^{2}\left( t\right) }V_{2}\left( .,t\right) \right\Vert _{B}<\infty 
\]%
or%
\[
\lim\limits_{R\rightarrow \infty }\left\Vert V\right\Vert _{L^{1}\left(
0,1;L^{\infty }\left( R^{n}/B_{R}\right) ;L\left( H\right) \right) }=0; 
\]

(3) $\alpha $, $\beta >0$, $\alpha \beta <2$ and $u\in C\left( \left[ 0,1%
\right] ;X\right) $ is a solution of the equation $\left( 2.6\right) $ and%
\[
\left\Vert e^{\frac{\left\vert x\right\vert ^{2}}{\beta ^{2}}}u\left(
.,0\right) \right\Vert _{X}<\infty ,\left\Vert e^{\frac{\left\vert
x\right\vert ^{2}}{\alpha ^{2}}}u\left( .,1\right) \right\Vert _{X}<\infty . 
\]

Then $u\left( x,t\right) \equiv 0.$

As a direct consequence of Theorem 1 we get the following Hardy's
uncertainty principle result for the non-linear equations (1.2).

\textbf{Theorem 2. }Let $u_{1},u_{2}\in C\left( \left[ 0,1\right]
;Y^{k}\right) $, $k\in Z^{+}$ be stronge solutions of the equation $\left(
1.2\right) $ with $k>\frac{n}{2}.$ Moreover, assume $F\in C^{k}\left( 
\mathbb{C}^{2},\mathbb{C}\right) $ and $F\left( 0\right) =\partial
_{u}F\left( 0\right) =\partial _{\bar{u}}F\left( 0\right) =0.$ If there are $%
\alpha $, $\beta >0$ with $\alpha \beta <2$ such that 
\[
e^{-\frac{\left\vert x\right\vert ^{2}}{\beta ^{2}}}\left( u_{1}\left(
.,0\right) -u_{2}\left( .,0\right) \right) \in X\text{, }e^{-\frac{%
\left\vert x\right\vert ^{2}}{\alpha ^{2}}}\left( u_{1}\left( .,1\right)
-u_{2}\left( .,1\right) \right) \in X 
\]%
then $u_{1}\equiv u_{2}.$

One of the results we get is the following one.

\textbf{Theorem 3. }Assume all conditions of Theorem 1 are satisfied. Then $%
\left\Vert e^{\left\vert x\right\vert ^{2}\mu ^{2}\left( t\right) }u\left(
.,t\right) \right\Vert _{X}^{\frac{1}{\mu \left( t\right) }}$ is
logarithmically convex in $[0,1]$ and there is $N=N\left( \alpha ,\beta
\right) $ such that 
\[
\left\Vert e^{\left\vert x\right\vert ^{2}\mu ^{2}\left( t\right) }u\left(
.,t\right) \right\Vert _{X}^{\frac{1}{\mu \left( t\right) }}\leq 
\]%
\[
e^{N\left( M_{1}+M_{2}+M_{1}^{2}+M_{2}^{2}\right) }\left\Vert e^{\frac{%
\left\vert x\right\vert ^{2}}{\beta ^{2}}}u\left( .,0\right) \right\Vert
_{X}^{\beta \left( 1-t\right) \mu \left( t\right) }\left\Vert e^{\frac{%
\left\vert x\right\vert ^{2}}{\alpha ^{2}}}u\left( .,1\right) \right\Vert
_{X}^{\alpha t\mu \left( t\right) }, 
\]%
when 
\[
M_{2}=\sup\limits_{t\in \left[ 0,1\right] }\left\Vert e^{\left\vert
x\right\vert ^{2}\mu ^{2}\left( t\right) }V_{2}\left( .,t\right) \right\Vert
_{B}e^{2\sup\limits_{t\in \left[ 0,1\right] }\left\Vert \func{Re}V_{2}\left(
.,t\right) \right\Vert _{B}}. 
\]%
Moreover,

\[
\sqrt{t\left( 1-t\right) }\left\Vert e^{\left\vert x\right\vert ^{2}\mu
^{2}\left( t\right) }\nabla u\right\Vert _{L^{2}\left( R^{n}\times \left[ 0,1%
\right] ;H\right) }\leq 
\]%
\[
e^{N\left( M_{1}+M_{2}+M_{1}^{2}+M_{2}^{2}\right) }\left[ \left\Vert e^{%
\frac{\left\vert x\right\vert ^{2}}{\beta ^{2}}}u\left( .,0\right)
\right\Vert _{X}+\left\Vert e^{\frac{\left\vert x\right\vert ^{2}}{\alpha
^{2}}}u\left( .,1\right) \right\Vert _{X}\right] . 
\]%
Here, we prove the following result for abstract parabolic equations with
variable coefficientes.

Consider the Cauchy problem for parabolic equations with variable operator
coefficients 
\begin{equation}
\partial _{t}u=\Delta u+Au+V\left( x,t\right) u,\text{ }x\in R^{n},\text{ }%
t\in \left[ 0,1\right] ,  \tag{2.7}
\end{equation}

\[
u\left( x,0\right) =f\left( x\right) , 
\]%
where $A$ is a linear operator and $V\left( x,t\right) $ is the given
potential operator function in $H$ $.$ By employing Theorem 1 we obtain

\textbf{Theorem 4. }Assume $A$ is a symmetric operator in $H$, $V\left(
x,t\right) \in L\left( H\right) $ for $\left( x,t\right) \in R^{n}\times %
\left[ 0,1\right] $\ and either, $V\left( x,t\right) =V_{1}\left( x\right)
+V_{2}\left( x,t\right) $, where $V_{1}\in L\left( H\right) $ and 
\[
\sup\limits_{t\in \left[ 0,1\right] }\left\Vert e^{\left\vert x\right\vert
^{2}\mu ^{2}\left( t\right) }V_{2}\left( .,t\right) \right\Vert _{B}<\infty 
\]
or 
\[
\lim\limits_{R\rightarrow \infty }\left\Vert V\right\Vert _{L^{1}\left(
0,1;L^{\infty }\left( R^{n}/O_{R}\right) ;L\left( H\right) \right) }=0. 
\]

Moreover, suppose $u\in L^{\infty }\left( 0,1;X\right) \cap L^{2}\left(
0,1;Y^{1}\right) $ is a solution of $\left( 2.9\right) $ and 
\[
\left\Vert f\right\Vert _{X}<\infty ,\left\Vert e^{\frac{\left\vert
x\right\vert ^{2}}{\delta ^{2}}}u\left( .,1\right) \right\Vert _{X}<\infty 
\]%
for some $\delta <1$. Then, $f\left( x\right) \equiv 0$ for $x\in R^{n}.$

\begin{center}
\textbf{3. Estimates for solutions}
\end{center}

\bigskip We need the following lemmas for proving the main results. Consider
the abstract Schr\"{o}dinger equation%
\begin{equation}
\partial _{t}u=\left( a+ib\right) \left[ \Delta u+Au+V\left( x,t\right)
u+F\left( x,t\right) \right] ,\text{ }x\in R^{n},\text{ }t\in \left[ 0,1%
\right] ,  \tag{3.0}
\end{equation}%
where $a$, $b$ are real numbers, $A$ is a linear operator, $V\left(
x,t\right) $ is a given potential operator function in $H$ and $F\left(
x,t\right) $ is a given $H$-valued function.

Let 
\[
\Phi \left( A,V\right) \upsilon =a\func{Re}\left( \left( A+V\right) \upsilon
,\upsilon \right) _{H}-b\func{Im}\left( \left( A+V\right) \upsilon ,\upsilon
\right) _{H}, 
\]%
\[
\text{ for }\upsilon =\upsilon \left( x,t\right) \in H\left( A\right) . 
\]

\textbf{Lemma 3.1. }Assume $a>0,$ $b\in \mathbb{R}$, $A$ is a symmetric
operator in $H.$ Moreover, there is a constant $C_{0}>0$ so that 
\begin{equation}
\left\vert \Phi \left( A,V\right) \upsilon \left( x,t\right) \right\vert
\leq C_{0}\mu \left( x,t\right) \left\Vert \upsilon \left( x,t\right)
\right\Vert _{H}^{2},  \tag{3.V}
\end{equation}%
for $x\in R^{n},$ $t\in \left[ 0,T\right] ,$ $\gamma \geq 0$, $T\in \left[
0,1\right] $ and $\upsilon \in H\left( A\right) $, where $\mu $ is a
positive function in $L^{1}\left( 0,T;L^{\infty }\left( R^{n}\right) \right) 
$.

Then the solution $u$ of $\left( 3.0\right) $ belonging to $L^{\infty
}\left( 0,1;X\right) \cap L^{2}\left( 0,1;Y^{1}\right) $ satisfies the
following estimate%
\[
e^{M_{T}}\left\Vert e^{\phi \left( .,T\right) }u\left( .,T\right)
\right\Vert _{X}\leq M_{T}\left\Vert e^{\gamma \left\vert x\right\vert
^{2}}u\left( .,0\right) \right\Vert _{X}+\sqrt{a^{2}+b^{2}}\left\Vert
e^{\phi \left( t\right) }F\right\Vert _{L^{1}\left( 0,T;X\right) }, 
\]%
where 
\[
\phi \left( x,t\right) =\frac{\gamma a\left\vert x\right\vert ^{2}}{%
a+4\gamma \left( a^{2}+b^{2}\right) t}\text{, }M_{T}=\left\Vert \mu
\right\Vert _{L^{1}\left( 0,T:L^{\infty }\left( R^{n}\right) \right) }. 
\]

\textbf{Proof. }Let $\upsilon =e^{\varphi }u$ where $\varphi $ is a
real-valued function to be chosen later. The function $\upsilon $ verifies%
\[
\partial _{t}\upsilon =S\upsilon +K\upsilon +\left( a+ib\right) \left[
\left( A+V\right) +e^{\varphi }F\right] \text{ in }R^{n}\times \left[ 0,1%
\right] \text{,} 
\]%
where $S$, $K$ are symmetric and skew-symmetric operators given by%
\[
S=a\left( \Delta +\left\vert \nabla \varphi \right\vert ^{2}\right)
-ib\left( 2\nabla \varphi .\nabla +\Delta \varphi \right) +\partial
_{t}\varphi , 
\]%
\[
K=ib\left( \Delta +\left\vert \nabla \varphi \right\vert ^{2}\right)
-a\left( 2\nabla \varphi .\nabla +\Delta \varphi \right) . 
\]%
By differentating inner product in $X$, we get 
\begin{equation}
\partial _{t}\left\Vert \upsilon \right\Vert _{X}^{2}=2\func{Re}\left(
S\upsilon ,\upsilon \right) _{X}+2\func{Re}\left( K\upsilon ,\upsilon
\right) _{X}+  \tag{3.1}
\end{equation}%
\[
2\func{Re}\left( \left( a+ib\right) e^{\upsilon }F,\upsilon \right) _{X}+2%
\func{Re}\left( \left( a+ib\right) \left( A+V\right) \upsilon ,\upsilon
\right) _{X},\text{ }t\geq 0. 
\]%
A formal integration by parts gives that%
\[
\func{Re}\left( S\upsilon ,\upsilon \right)
_{X}=-a\dint\limits_{R^{n}}\left\vert \nabla \upsilon \right\vert
_{H}^{2}dx+\dint\limits_{R^{n}}\left( a\left\vert \nabla \varphi \right\vert
^{2}+\partial _{t}\varphi \right) \left\Vert \upsilon \right\Vert
_{H}^{2}dx+ 
\]%
\[
b\func{Im}\dint\limits_{R^{n}}\left( 2\nabla \varphi .\Delta \upsilon
,\upsilon \right) _{H}dx, 
\]%
\begin{equation}
\func{Re}\left( K\upsilon ,\upsilon \right) _{X}=\left( -2a\nabla \varphi
.\left( \nabla \upsilon ,\upsilon \right) \right) _{X}-a\func{Re}%
\dint\limits_{R^{n}}\left( 2\nabla \varphi .\Delta \upsilon ,\upsilon
\right) _{H}dx,  \tag{3.2}
\end{equation}%
\[
\func{Re}\left( \left( a+ib\right) \left( A+V\right) \upsilon ,\upsilon
\right) _{X}=a\func{Re}\dint\limits_{R^{n}}\left( \left( A+V\right) \upsilon
,\upsilon \right) _{H}dx 
\]

\[
-b\func{Im}\dint\limits_{R^{n}}\left( \left( A+V\right) \upsilon ,\upsilon
\right) _{H}dx=\dint\limits_{R^{n}}\Phi \left( A,V\right) \upsilon dx, 
\]

\[
\func{Re}\left( \left( a+ib\right) e^{\varphi }F,\upsilon \right) _{X}=a%
\func{Re}\dint\limits_{R^{n}}\left( e^{\varphi }F,\upsilon \right)
_{H}dx=ae^{\varphi }\func{Re}\left( F,\upsilon \right) _{X}. 
\]

By using the Cauchy-Schwarz's inequality, by condition $\left( 3.V\right) $,
in view of $\left( 3.1\right) $ and $\left( 3.2\right) $ we obtain%
\[
\partial _{t}\left\Vert \upsilon \right\Vert ^{2}\leq ae^{\varphi
}\left\Vert F\left( t,.\right) \right\Vert _{X}\left\Vert \upsilon
\right\Vert _{X}+C_{0}\left\Vert \mu \left( .,t\right) \right\Vert
_{L^{\infty }\left( R^{n}\right) }\left\Vert \upsilon \right\Vert _{X}^{2}, 
\]%
where $a,$ $b$ and $\varphi $ are such that%
\begin{equation}
\left( a-b\right) \Delta \varphi \leq 0,\text{ }\left( a+\frac{b^{2}}{a}%
\right) \left\vert \nabla \varphi \right\vert ^{2}+\partial _{t}\upsilon
\leq 0\text{ in }R_{+}^{n+1}.  \tag{3.3}
\end{equation}%
The remainig part of the proof is obtained by reasoning as in $\left[ \text{%
7, Lemma 1}\right] .$

When $\varphi \left( x,t\right) =q\left( t\right) \psi \left( x\right) $, it
suffices that%
\begin{equation}
\left( a+\frac{b^{2}}{a}\right) q^{2}\left( t\right) \left\vert \nabla \psi
\right\vert ^{2}+q^{\prime }\left( t\right) \psi \left( x\right) \leq 0. 
\tag{3.4}
\end{equation}%
If we put $\psi \left( x\right) =\left\vert x\right\vert ^{2}$ then $\left(
3.4\right) $ holds, when 
\begin{equation}
q^{\prime }\left( t\right) =-4\left( a+\frac{b^{2}}{a}\right) q^{2}\left(
t\right) ,\text{ }q\left( 0\right) =\gamma ,\text{ }\gamma \geq 0.  \tag{3.5}
\end{equation}

Let 
\[
\psi _{R}\left( x\right) =\left\{ 
\begin{array}{c}
\left\vert x\right\vert ^{2}\text{, }\left\vert x\right\vert <R \\ 
\infty \text{, \ \ \ \ }\left\vert x\right\vert >R\text{\ \ \ \ \ \ \ \ }%
\end{array}%
.\right. 
\]%
Regularize $\psi _{R}$ with a radial mollifier $\theta _{\rho }$ and set 
\[
\varphi _{\rho ,R}\left( x,t\right) =q\left( t\right) \theta _{\rho }\ast
\psi _{R}\left( x\right) ,\text{ }\upsilon _{\rho ,R}\left( x,t\right)
=e^{\varphi _{\rho ,R}}u,\text{ } 
\]%
where $q\left( t\right) =\gamma a\left[ a+4\gamma \left( a^{2}+b^{2}\right) t%
\right] ^{-1}$ is the solution to $\left( 3.5\right) $. Because the right
hand side of $(3.2)$ only involves the first derivatives of $\varphi $, $%
\psi _{R}$ is Lipschitz and bounded at infinity, 
\[
\theta _{\rho }\ast \psi _{R}\left( x\right) \leq \theta _{\rho }\ast
\left\vert x\right\vert ^{2}=C\left( n\right) \rho ^{2} 
\]%
and $(3.3)$ holds uniformly in $\rho $ and $R$, when $\varphi $ is replaced
by $\varphi _{\rho ,R}$. Hence, it follows that the estimate%
\[
e^{M_{T}}\left\Vert e^{\phi \left( T\right) }u\left( T\right) \right\Vert
_{X}\leq M_{T}\left\Vert e^{\gamma \left\vert x\right\vert ^{2}}u\left(
0\right) \right\Vert _{X}+\sqrt{a^{2}+b^{2}}\left\Vert e^{\varphi _{\rho
,R}}F\right\Vert _{L^{1}\left( 0,T;X\right) } 
\]%
holds uniformly in $\rho $ and $R.$ The assertion is obtained after letting $%
\rho $ tend to zero and $R$ to infinity.

\textbf{Remark 3.1. }It should be noted that if $H=\mathbb{C}$, $A=0$ and $%
V\left( x,t\right) $ is a complex valued function, then the abstract
condition $\left( 3.V\right) $ can be replised by 
\[
M_{T}=\left\Vert a\func{Re}V-b\func{Im}V\right\Vert _{L^{1}\left(
0,T;L^{\infty }\left( R^{n}\right) \right) }<\infty . 
\]%
Moreover, if $A$ and $V\left( x,t\right) $ for $x\in R^{n},$ $t\in \left[ 0,T%
\right] $ are bounded operators in $H$, then by using Cauchy-Schwarz's
inequality the assumption $\left( 3.V\right) $ replaced as 
\[
\left\vert \Phi \left( A,V\right) \upsilon \right\vert \leq \sqrt{a^{2}+b^{2}%
}\left\Vert \left( A+V\right) \upsilon \right\Vert _{H}\left\Vert \upsilon
\right\Vert _{H}\leq \sqrt{a^{2}+b^{2}}\left\Vert A+V\right\Vert _{L\left(
H\right) }\left\Vert \upsilon \right\Vert ^{2}. 
\]

Let 
\[
\text{ }Q\left( t\right) =\left( f,f\right) _{X}\text{, }D\left( t\right)
=\left( Sf,f\right) _{X},\text{ }N\left( t\right) =D\left( t\right)
Q^{-1}\left( t\right) , 
\]
\[
\text{ }\partial _{t}S=S_{t}\text{ and }\left\vert \nabla \upsilon
\right\vert _{H}^{2}=\dsum\limits_{k=1}^{n}\left\Vert \frac{\partial
\upsilon }{\partial x_{k}}\right\Vert _{H}^{2}. 
\]

\textbf{Lemma 3.2. }Assume\ $S=S\left( t\right) $ is a symmetric, $K=K\left(
t\right) $ is a skew-symmetric operators in $H$, $G\left( x,t\right) $ is a
positive funtion and $f(x,t)$ is a reasonable function. Then, 
\[
Q^{^{\prime \prime }}\left( t\right) =2\partial _{t}\func{Re}\left( \partial
_{t}f-Sf-Kf,f\right) _{X}+2\left( S_{t}f+\left[ S,K\right] f,f\right) _{X}+ 
\]

\begin{equation}
\left\Vert \partial _{t}f-Sf+Kf\right\Vert _{X}^{2}-\left\Vert \partial
_{t}f-Sf-Kf\right\Vert _{X}^{2}  \tag{3.6}
\end{equation}%
and 
\[
\partial _{t}N\left( t\right) \geq Q^{-1}\left( t\right) \left[ \left(
S_{t}f+\left[ S,K\right] f,f\right) _{X}-\frac{1}{2}\left\Vert \partial
_{t}f-Sf-Kf\right\Vert _{X}^{2}\right] . 
\]

\bigskip Moreover, if 
\[
\left\Vert \partial _{t}f-Sf-Kf\right\Vert _{H}\leq M_{1}\left\Vert
f\right\Vert _{H}+G\left( x,t\right) ,\text{ }S_{t}+\left[ S,K\right] \geq
-M_{0}\text{ } 
\]%
for $x\in R^{N}$, $t\in \left[ 0,1\right] $ and%
\[
M_{2}=\sup\limits_{t\in \left[ 0,1\right] }\left\Vert G\left( .,t\right)
\right\Vert _{L^{2}\left( R^{n}\right) }\left( \left\Vert f\left( .,t\right)
\right\Vert _{X}\right) ^{-1}<\infty . 
\]

Then $Q\left( t\right) $ is logarithmically convex in $[0,1]$ and there is a
constant $M$ such that 
\[
Q\left( t\right) \leq e^{M\left(
M_{0}+M_{1}+M_{2}+M_{1}^{2}+M_{2}^{2}\right) }Q^{1-t}\left( 0\right)
Q^{t}\left( 1\right) \text{, }0\leq t\leq 1. 
\]

\textbf{Proof.} The lemma is verifying in a similar way as in $\left[ \text{%
7, Lemma 2}\right] $ by replacing the inner product and norm of $L^{2}\left(
R^{n}\right) $ with inner product and norm of the space $L^{2}\left(
R^{n};H\right) .$

\bigskip \textbf{Lemma 3.3. }Assume $a$, $\gamma >0,$ $b\in \mathbb{R}$, $A$
is a symmetric operator in $H.$ Let, 
\[
\left\vert \Phi \left( A,V\right) u\left( x,t\right) \right\vert \leq
C_{0}\mu \left( x,t\right) \left\Vert u\left( x,t\right) \right\Vert
_{H}^{2} 
\]%
for $x\in R^{n},$ $t\in \left[ 0,1\right] $ and $u\in H\left( A\right) $,
where $\mu $ is a positive function in $L^{1}\left( 0,T;L^{\infty }\left(
R^{n}\right) \right) $. Moreover, suppose%
\[
\sup\limits_{t\in \left[ 0,1\right] }\left\Vert V\left( .,t\right)
\right\Vert _{B}\leq M_{1}\text{, }\left\Vert e^{\gamma \left\vert
x\right\vert ^{2}}u\left( .,0\right) \right\Vert _{X}<\infty ,\text{ }%
\left\Vert e^{\gamma \left\vert x\right\vert ^{2}}u\left( .,1\right)
\right\Vert _{X}<\infty \text{ } 
\]%
and 
\[
M_{2}=\sup\limits_{t\in \left[ 0,1\right] }\frac{\left\Vert e^{\gamma
\left\vert x\right\vert ^{2}}F\left( .,t\right) \right\Vert _{X}}{\left\Vert
u\right\Vert _{X}}<\infty . 
\]

Then, for solution $u\in L^{\infty }\left( 0,1;X\right) \cap L^{2}\left(
0,1;Y^{1}\right) $ of the equation $\left( 3.0\right) $, $e^{\gamma
\left\vert x\right\vert ^{2}}u\left( .,t\right) $ is logarithmically convex
in $[0,1]$ and there is a constant $N$ such that%
\begin{equation}
\left\Vert e^{\gamma \left\vert x\right\vert ^{2}}u\left( .,t\right)
\right\Vert _{X}\leq e^{NM\left( a,b\right) }\left\Vert e^{\gamma \left\vert
x\right\vert ^{2}}u\left( .,0\right) \right\Vert _{X}^{1-t}\left\Vert
e^{\gamma \left\vert x\right\vert ^{2}}u\left( .,1\right) \right\Vert
_{X}^{t}  \tag{3.7}
\end{equation}%
where 
\[
M\left( a,b\right) =\left( a^{2}+b^{2}\right) \left( \gamma
M_{1}^{2}+M_{2}^{2}\right) +\sqrt{a^{2}+b^{2}}\left( M_{1}+M_{2}\right) 
\]%
when $0\leq t\leq 1.$

\textbf{Proof. }Let $f=e^{\gamma \varphi }u,$ where $\varphi $ is a
real-valued function to be chosen. The function $f\left( x\right) $ verifies%
\begin{equation}
\partial _{t}f=Sf+Kf+\left( a+ib\right) \left( Vf+e^{\gamma \varphi
}F\right) \text{ in }R^{n}\times \left[ 0,1\right] \text{,}  \tag{3.8}
\end{equation}%
where $S$, $K$ are symmetric and skew-symmetric operator, respectively given
by%
\begin{equation}
S=a\left( \Delta +A+\gamma ^{2}\left\vert \nabla \varphi \right\vert
^{2}\right) -ib\gamma \left( 2\nabla \varphi .\nabla +\Delta \varphi \right)
+\gamma \partial _{t}\varphi ,  \tag{3.9}
\end{equation}%
\[
K=ib\left( \Delta +A+\gamma ^{2}\left\vert \nabla \varphi \right\vert
^{2}\right) -a\gamma \left( 2\nabla \varphi .\nabla +\Delta \varphi \right)
. 
\]%
A calculation shows that, 
\[
S_{t}+\left[ S,K\right] =\gamma \partial _{t}^{2}\varphi +4\gamma
^{2}a\nabla \varphi .\nabla \partial _{t}\varphi -2ib\gamma \left( 2\nabla
\partial _{t}\varphi .\nabla +\Delta \partial _{t}\varphi \right) - 
\]%
\begin{equation}
\gamma \left( a^{2}+b^{2}\right) \left[ 4\nabla .\left( D^{2}\varphi \nabla
\right) -4\gamma ^{2}D^{2}\varphi \nabla \varphi +\Delta ^{2}\varphi \right]
.  \tag{3.10}
\end{equation}%
If we put $\varphi =\left\vert x\right\vert ^{2}$, then $\left( 3.10\right) $
reduce the following%
\[
S_{t}+\left[ S,K\right] =-\gamma \left( a^{2}+b^{2}\right) \left[ 8\Delta
-32\gamma ^{2}\left\vert x\right\vert ^{2}\right] . 
\]%
Moreover, 
\begin{equation}
\left( S_{t}f+\left[ S,K\right] f,f\right) =\gamma \left( a^{2}+b^{2}\right)
\dint\limits_{R^{n}}\left( 8\left\vert \nabla f\right\vert _{H}^{2}+32\gamma
^{2}\left\vert x\right\vert ^{2}\left\Vert f\right\Vert _{H}^{2}\right) dx. 
\tag{3.11}
\end{equation}%
This identity, the condition on $V$ and $\left( 3.8\right) $ imply that%
\begin{equation}
\left\Vert \partial _{t}f-Sf-Kf\right\Vert _{X}\leq \sqrt{a^{2}+b^{2}}\left(
M_{1}\left\Vert f\right\Vert _{X}+e^{\gamma \varphi }\left\Vert F\right\Vert
_{X}\right) \text{, }S_{t}+\left[ S,K\right] \geq 0.  \tag{3.12}
\end{equation}%
If we knew that the quantities and calculations involved in the proof of
Lemma 3.2 (similar as $\left[ \text{7, Lemma 2}\right] $) were finite and
correct, when $f=e^{\gamma \left\vert x\right\vert ^{2}}u$ we would have the
logarithmic convexity of $Q\left( t\right) =\left\Vert e^{\gamma \left\vert
x\right\vert ^{2}}u\left( .,t\right) \right\Vert _{X}$ \ and the estimate $%
\left( 3.7\right) $ from Lemma 3.2. But this fact is verifying by reasonong
as in $\left[ \text{7, Lemma 3}\right] .$

Let 
\[
\eta =\sqrt{t\left( 1-t\right) }e^{\gamma \left\vert x\right\vert ^{2}},%
\text{ }Z=L^{2}\left( \left[ 0,1\right] \times R^{n};H\right) . 
\]

\textbf{Lemma 3.4. }Assume that $a$, $b$, $u$, $A$ and $V$ are as in Lemma
3.3 and $\gamma >0$. Then,

\ \ \ \ \ \ \ \ \ \ \ \ \ \ \ \ \ \ \ \ \ \ \ \ \ \ \ \ \ \ \ \ \ \ \ \ \ \
\ \ \ \ \ \ \ \ \ \ \ \ \ \ \ \ \ \ \ \ \ \ \ \ \ \ \ \ \ \ \ \ \ \ \ \ \ \
\ \ \ \ \ \ \ \ \ \ \ \ \ \ \ 
\[
\left\Vert \eta \nabla u\right\Vert _{Z}+\left\Vert \eta \left\vert
x\right\vert u\right\Vert _{Z}\leq N\left[ \left( 1+M_{1}\right) \right] %
\left[ \sup\limits_{t\in \left[ 0,1\right] }\left\Vert e^{\gamma \left\vert
x\right\vert ^{2}}u\left( .,t\right) \right\Vert _{X}+\sup\limits_{t\in %
\left[ 0,1\right] }\left\Vert e^{\gamma \left\vert x\right\vert ^{2}}F\left(
.,t\right) \right\Vert _{Z}\right] \text{,} 
\]%
where $N$ is bounded number, when $\gamma \ $and $a^{2}+b^{2}$ are bounded
below.

\textbf{Proof. }The integration by parts shows that 
\[
\dint\limits_{R^{n}}\left( \left\vert \nabla f\right\vert _{H}^{2}+4\gamma
^{2}\left\vert x\right\vert ^{2}\left\Vert f\right\Vert _{H}\right)
dx=\dint\limits_{R^{n}}\left[ e^{2\gamma \left\vert x\right\vert ^{2}}\left(
\left\vert \nabla u\right\vert _{H}^{2}-2n\gamma \right) \left\Vert
u\right\Vert _{H}^{2}\right] dx, 
\]%
when $f=e^{\gamma \left\vert x\right\vert ^{2}}u$, while integration by
parts, the Cauchy-Schwarz's inequality and the identity, $n=\nabla $%
\textperiodcentered $x$, give that 
\[
\dint\limits_{R^{n}}\left( \left\vert \nabla f\right\vert _{H}^{2}+4\gamma
^{2}\left\vert x\right\vert ^{2}\left\Vert f\right\Vert _{H}\right) dx\geq
2\gamma n\left\Vert f\right\Vert _{X}^{2}\text{.} 
\]

The sum of the last two formulae gives the inequality%
\begin{equation}
2\dint\limits_{R^{n}}\left( \left\vert \nabla f\right\vert _{H}^{2}+4\gamma
^{2}\left\vert x\right\vert ^{2}\left\Vert f\right\Vert _{H}\right) dx\geq
\dint\limits_{R^{n}}e^{\gamma \left\vert x\right\vert ^{2}}\left\vert \nabla
f\right\vert _{H}^{2}dx\text{.}  \tag{3.13}
\end{equation}%
Integration over $[0,1]$ of $t(1-t)$ times the formula (3.6) for $%
Q^{^{\prime \prime }}\left( t\right) $ and integration by parts, shows that 
\begin{equation}
2\dint\limits_{0}^{1}t(1-t)\left( S_{t}f+\left[ S,K\right] f,f\right)
_{X}dt+\dint\limits_{0}^{1}Q\left( t\right) dt\leq Q\left( 1\right) +Q\left(
0\right) +  \tag{3.14}
\end{equation}

\[
2\dint\limits_{0}^{1}1-2t)\func{Re}\left( \partial _{t}f-Sf-Kf,f\right)
_{X}dx+\dint\limits_{0}^{1}t(1-t)\left\Vert \partial _{t}f-Sf-Kf\right\Vert
_{X}^{2}dt. 
\]

Assuming again that the last two calculations are justified for $f=e^{\gamma
\left\vert x\right\vert ^{2}}.$ Then $\left( 2.11\right) -\left( 2.14\right) 
$ implay the assertion.

\begin{center}
\textbf{4. Appell transformation in abstract functon spaces}
\end{center}

\bigskip Let 
\[
\text{ }\eta \left( x,t\right) =\frac{\left( \alpha -\beta \right)
\left\vert x^{2}\right\vert }{4\left( a+ib\right) \alpha \left( 1-t\right)
+\beta t},\text{ }\nu \left( s\right) =\left[ \gamma \alpha \beta \mu
^{2}\left( s\right) +\frac{\left( \alpha -\beta \right) a}{4\left(
a^{2}+b^{2}\right) }\mu \left( s\right) \right] . 
\]

\bigskip \textbf{Lemma 4.1. }Assume $A$ and $V$ are as in Lemma 3.3 and $%
u=u\left( x,s\right) $ is a solution of the equation 
\[
\partial _{s}u=\left( a+ib\right) \left[ \Delta u+Au+V\left( y,s\right)
u+F\left( y,s\right) \right] ,\text{ }y\in R^{n},\text{ }s\in \left[ 0,1%
\right] . 
\]%
Let $a+ib\neq 0$, $\gamma \in \mathbb{R}$ and $\alpha $, $\beta \in \mathbb{R%
}_{+}$. Set 
\begin{equation}
\tilde{u}\left( x,t\right) =\left( \sqrt{\alpha \beta }\mu \left( t\right)
\right) ^{\frac{n}{2}}u\left( \sqrt{\alpha \beta }x\mu \left( t\right)
,\beta t\mu \left( t\right) \right) e^{\eta }.  \tag{4.1}
\end{equation}

\bigskip Then, $\tilde{u}\left( x,t\right) $\ verifies the equation 
\[
\partial _{t}\tilde{u}=\left( a+ib\right) \left[ \Delta \tilde{u}+A\tilde{u}+%
\tilde{V}\left( x,t\right) u+\tilde{F}\left( x,t\right) \right] ,\text{ }%
x\in R^{n},\text{ }t\in \left[ 0,1\right] 
\]%
with 
\[
\tilde{V}\left( x,t\right) =\alpha \beta \mu ^{2}\left( t\right) V\left( 
\sqrt{\alpha \beta }x\mu \left( t\right) ,\beta t\mu \left( t\right) \right)
, 
\]

\[
\text{ }\tilde{F}\left( x,t\right) =\left( \sqrt{\alpha \beta }\mu \left(
t\right) \right) ^{\frac{n}{2}+2}\left( \sqrt{\alpha \beta }x\mu \left(
t\right) ,\beta t\mu \left( t\right) \right) . 
\]%
Moreover, 
\[
\left\Vert e^{\gamma \left\vert x\right\vert ^{2}}\tilde{F}\left( .,t\right)
\right\Vert _{X}=\alpha \beta \mu ^{2}\left( t\right) e^{\nu \left\vert
y\right\vert ^{2}}\left\Vert F\left( s\right) \right\Vert _{X}\text{ and }%
\left\Vert e^{\gamma \left\vert x\right\vert ^{2}}\tilde{u}\left( .,t\right)
\right\Vert _{X}=e^{\nu \left\vert y\right\vert ^{2}}\left\Vert u\left(
s\right) \right\Vert _{X} 
\]%
when $s=\mu \left( t\right) $ and $\gamma \in \mathbb{R}$.

\textbf{Proof. }If $u$ is a solution of the equation%
\begin{equation}
\partial _{s}u=\left( a+ib\right) \left[ \Delta u+Au+Q\left( y,s\right) %
\right] ,\text{ }y\in R^{n},\text{ }s\in \left[ 0,1\right]  \tag{4.2}
\end{equation}%
then, the function $u_{1}\left( x,t\right) =u\left( \sqrt{r}x,rt+\tau
\right) $ verifies%
\[
\partial _{t}u_{1}=\left( a+ib\right) \left[ \Delta u_{1}+Au_{1}+rQ\left( 
\sqrt{r}x,rt+\tau \right) \right] ,\text{ }y\in R^{n},\text{ }s\in \left[ 0,1%
\right] 
\]%
and $u_{2}\left( x,t\right) =t^{-\frac{n}{2}}u\left( \frac{x}{t},\frac{1}{t}%
\right) e^{\frac{\left\vert x\right\vert ^{2}}{4\left( a+ib\right) t}}$ is a
\ solution to%
\[
\partial _{t}u_{2}=-\left( a+ib\right) \left[ \Delta u_{2}+Au+t^{-\left( 2+%
\frac{n}{2}\right) }Q\left( \frac{x}{t},\frac{1}{t}\right) e^{\frac{%
\left\vert x\right\vert ^{2}}{4\left( a+ib\right) t}}\right] ,\text{ }y\in
R^{n},\text{ }s\in \left[ 0,1\right] . 
\]%
These two facts and the sequel of changes of variables below verifies the
Lemma, when $\alpha >\beta ,$ i.e.%
\[
u\left( \sqrt{\frac{\alpha \beta }{\alpha -\beta }}x,\frac{\alpha \beta }{%
\alpha -\beta }t-\frac{\beta }{\alpha -\beta }\right) 
\]%
\ \ is a solution to the same non-homogeneous equation but with right-hand
side 
\[
\frac{\alpha \beta }{\alpha -\beta }Q\left( \sqrt{\frac{\alpha \beta }{%
\alpha -\beta }}x,\frac{\alpha \beta }{\alpha -\beta }t-\frac{\beta }{\alpha
-\beta }\right) . 
\]

The function, 
\[
\frac{1}{\left( \alpha -t\right) ^{\frac{n}{2}}}u\left( \frac{\sqrt{\alpha
\beta }x}{\sqrt{\alpha -\beta }\left( \alpha -t\right) },\frac{\alpha \beta 
}{\left( \alpha -\beta \right) \left( \alpha -t\right) }-\frac{\beta }{%
\alpha -\beta }\right) e^{\frac{\left\vert x\right\vert ^{2}}{4\left(
a+ib\right) \left( \alpha -t\right) }} 
\]%
verifies $(4.2)$ with right-hand side%
\[
\frac{\alpha \beta }{\left( \alpha -\beta \right) \left( \alpha -t\right) ^{%
\frac{n}{2}+2}}Q\left( \frac{\sqrt{\alpha \beta }x}{\sqrt{\alpha -\beta }%
\left( \alpha -t\right) },\frac{\alpha \beta }{\left( \alpha -\beta \right)
\left( \alpha -t\right) }-\frac{\beta }{\alpha -\beta }\right) e^{\frac{%
\left\vert x\right\vert ^{2}}{4\left( a+ib\right) \left( \alpha -t\right) }%
}. 
\]%
Replacing $\left( x,t\right) $ by $\left( \sqrt{\alpha -\beta }x,\left(
\alpha -\beta \right) t\right) $ we get that%
\begin{equation}
\mu ^{\frac{n}{2}}\left( t\right) u\left( \sqrt{\alpha \beta }\mu \left(
t\right) x,\frac{\alpha \beta \mu \left( t\right) }{\left( \alpha -\beta
\right) }-\frac{\beta }{\alpha -\beta }\right) e^{^{\left( \alpha -\beta
\right) }\frac{\left\vert x\right\vert ^{2}\mu \left( t\right) }{4\left(
a+ib\right) }}  \tag{4.3}
\end{equation}%
is a solution of $\left( 4.2\right) $ but with right-hand%
\begin{equation}
\mu ^{\frac{n}{2}+2}\left( t\right) Q\left( \sqrt{\alpha \beta }\mu \left(
t\right) x,\frac{\alpha \beta \mu \left( t\right) }{\left( \alpha -\beta
\right) }-\frac{\beta }{\alpha -\beta }\right) e^{^{\left( \alpha -\beta
\right) }\frac{\left\vert x\right\vert ^{2}\mu \left( t\right) }{4\left(
a+ib\right) }}.  \tag{4.4}
\end{equation}%
Finally, observe that%
\[
s=\beta t\mu \left( t\right) =\frac{\alpha \beta \mu \left( t\right) }{%
\left( \alpha -\beta \right) }-\frac{\beta }{\alpha -\beta } 
\]%
and multiply $\left( 4.3\right) $ and $\left( 4.4\right) $ we obtain the
assertion for $\alpha >\beta .$ The case $\beta >\alpha $ follows by
reversing by changes of variables, $s^{\prime }=1-s$ and $t^{\prime }=1-t.$

\begin{center}
\bigskip \textbf{\ 5. Variable coefficients. Proof of Theorem 3}
\end{center}

We are ready to prove Theorem 3.\textbf{\ }Let 
\[
B=L^{1}\left( 0,1;L^{\infty }\left( R^{n};L\left( H\right) \right) \right) ,%
\text{ }B\left( R\right) =L^{1}\left( 0,1;L^{\infty }\left(
R^{n}/O_{R};L\left( H\right) \right) \right) . 
\]

\textbf{Proof of Theorem 3. }We may assume that $\alpha \neq \beta $. The
case $\alpha =\beta $ follows from the latter by replacing $\beta $ by $%
\beta +\delta ,$ $\delta >0$, and letting $\delta $ tend to zero. We may
also assume that $\alpha <\beta $. Otherwise, replace $u$ by $\bar{u}(1-t)$.
Assume $a>0.$ Set $W=\Delta +A+V_{1}$ and let $Uu_{0}=e^{t\left( a+ib\right)
W}u_{0}$ denote $C\left( \left[ 0,1\right] ;X\right) $ solution to the
problem

\begin{equation}
\partial _{t}u=\left( a+ib\right) \left[ \Delta u+Au+V_{1}\left( x\right) u%
\right] ,\text{ }x\in R^{n},\text{ }t\in \left[ 0,1\right] ,  \tag{5.1}
\end{equation}%
\[
u\left( x,0\right) =u_{0}\left( x\right) . 
\]

By virtue of Duhamel principle there is a solution of 
\[
\partial _{t}u=\left( a+ib\right) \Delta u+Au+V\left( x,t\right) u,\text{ }%
x\in R^{n},\text{ }t\in \left[ 0,1\right] , 
\]%
\[
u\left( x,0\right) =u_{0}\left( x\right) . 
\]
expressing as 
\begin{equation}
u\left( x,t\right) =e^{t\left( a+ib\right)
W}u_{0}+i\dint\limits_{0}^{t}e^{-i\left( t-s\right) W}V_{2}\left( x,s\right)
u\left( x,s\right) ds\text{ }  \tag{5.2}
\end{equation}%
\[
\text{for }x\in R^{n},\text{ }s\in \left[ 0,1\right] . 
\]

For $0\leq \varepsilon \leq 1,$ set 
\begin{equation}
F_{\varepsilon }\left( x,t\right) =\frac{i}{\varepsilon +i}e^{\varepsilon
tW}V_{2}\left( x,t\right) u\left( x,t\right)  \tag{5.3}
\end{equation}

and 
\begin{equation}
u_{\varepsilon }\left( x,t\right) =e^{\left( \varepsilon +i\right)
tW}u_{0}+\left( \varepsilon +i\right) \dint\limits_{0}^{t}e^{\left(
\varepsilon +i\right) \left( t-s\right) W}F_{\varepsilon }\left( x,s\right)
u\left( x,s\right) ds.\text{ }  \tag{5.4}
\end{equation}

Then, $u_{\varepsilon }\left( x,t\right) \in L^{\infty }\left( 0,1;X\right)
\cap L^{2}\left( R^{n};Y^{1}\right) $ and satisfies 
\[
\partial _{t}u_{\varepsilon }=\left( \varepsilon +i\right) \left(
Wu+F_{\varepsilon }\right) \text{ in }R^{n}\times \left[ 0,1\right] , 
\]%
\[
u_{\varepsilon }\left( .,0\right) =u_{0}\left( .\right) . 
\]

The identities 
\[
e^{\left( z_{1}+z_{2}\right) W}=e^{z_{1}W}e^{z_{2}W}\text{, when }\func{Re}%
z_{1}\text{, }\func{Re}z_{2}\geq 0, 
\]

$\left( 5.2\right) $, $\left( 5.3\right) $ and $\left( 5.4\right) $ shows
that 
\begin{equation}
u_{\varepsilon }\left( x,t\right) =e^{\varepsilon tW}u\left( x,t\right) 
\text{, for }t\in \left[ 0,1\right] .  \tag{5.5}
\end{equation}%
In particular, the equality $u_{\varepsilon }\left( x,1\right)
=e^{\varepsilon W}u\left( x,1\right) $ and Lemma 3.1 with $a+ib=\varepsilon
, $ $\gamma =\frac{1}{\beta }$, $F\equiv 0$ and the fact that $%
u_{\varepsilon }(0)=u\left( 0\right) $ imply that%
\[
\left\Vert e^{\frac{\left\vert x\right\vert ^{2}}{\beta ^{2}+4\varepsilon }%
}u_{\varepsilon }\left( .,1\right) \right\Vert _{X}\leq e^{\varepsilon
\left\Vert V_{1}\right\Vert _{B}}\left\Vert e^{\frac{\left\vert x\right\vert
^{2}}{\beta ^{2}}}u\left( .,1\right) \right\Vert _{X}\text{, }\left\Vert e^{%
\frac{\left\vert x\right\vert ^{2}}{^{\alpha ^{2}}}}u_{\varepsilon }\left(
.,0\right) \right\Vert _{X}=\left\Vert e^{\frac{\left\vert x\right\vert ^{2}%
}{^{\alpha ^{2}}}}u\left( .,0\right) \right\Vert _{X}\text{.} 
\]%
A second application of Lemma 3.1 with $a+ib$ $=$ $\varepsilon $, $F\equiv 0$%
, the value of $\gamma =\mu ^{2}\left( t\right) $ and $\left( 5.2\right) $
show that%
\[
\left\Vert \varepsilon ^{\frac{\left\vert x\right\vert ^{2}}{\mu ^{2}\left(
t\right) +4\varepsilon t}}F_{\varepsilon }\left( .,t\right) \right\Vert
_{X}\leq e^{\varepsilon \left\Vert V_{1}\right\Vert _{B}}\left\Vert
\varepsilon ^{\frac{\left\vert x\right\vert ^{2}}{\mu ^{2}\left( t\right) }%
}V_{2}\left( .,t\right) \right\Vert _{B}\left\Vert u\left( .,t\right)
\right\Vert _{X}\text{, }t\in \left[ 0,1\right] . 
\]%
Setting, $\alpha _{\varepsilon }=\alpha +2\varepsilon $ and $\beta
_{\varepsilon }=\beta +2\varepsilon $, the last three inequalities give that%
\begin{equation}
\left\Vert e^{\frac{\left\vert x\right\vert ^{2}}{^{\beta _{\varepsilon
}^{2}}}}u_{\varepsilon }\left( .,1\right) \right\Vert _{X}\leq
e^{\varepsilon \left\Vert V_{1}\right\Vert _{B}}\left\Vert e^{\frac{%
\left\vert x\right\vert ^{2}}{^{\beta 2}}}u\left( .,1\right) \right\Vert
_{X},  \tag{5.6}
\end{equation}%
\[
\left\Vert e^{\frac{\left\vert x\right\vert ^{2}}{^{\alpha _{\varepsilon
}^{2}}}}u_{\varepsilon }\left( .,0\right) \right\Vert _{X}\leq
e^{\varepsilon \left\Vert V_{1}\right\Vert _{B}}\left\Vert e^{\frac{%
\left\vert x\right\vert ^{2}}{^{\alpha ^{2}}}}u\left( .,0\right) \right\Vert
_{X}, 
\]%
\begin{equation}
\left\Vert \varepsilon ^{\frac{\left\vert x\right\vert ^{2}}{\mu ^{2}\left(
t\right) }}F_{\varepsilon }\left( .,t\right) \right\Vert _{X}\leq
e^{\varepsilon \left\Vert V_{1}\right\Vert _{B}}\left\Vert \varepsilon ^{%
\frac{\left\vert x\right\vert ^{2}}{\mu ^{2}\left( t\right) }}V_{2}\left(
.,t\right) \right\Vert _{B}\left\Vert u\left( .,t\right) \right\Vert _{X}%
\text{, }t\in \left[ 0,1\right] .  \tag{5.7}
\end{equation}%
A third application of Lemma 3.1 with $a+ib=b,$ $F\equiv 0,$ $\gamma =0$,
and $(5.2),(5.5)$ implies that%
\begin{equation}
\left\Vert F_{\varepsilon }\left( .,t\right) \right\Vert _{X}\leq
e^{\varepsilon \left\Vert V_{1}\right\Vert _{B}}\left\Vert V_{2}\left(
.,t\right) \right\Vert _{B}\left\Vert u\left( .,t\right) \right\Vert _{X}%
\text{, }  \tag{5.8}
\end{equation}%
\[
\left\Vert u_{\varepsilon }\left( .,t\right) \right\Vert _{X}\leq
e^{\varepsilon \left\Vert V_{1}\right\Vert _{B}}\left\Vert u\left(
.,t\right) \right\Vert _{X}\text{, }t\in \left[ 0,1\right] . 
\]%
Set $\gamma _{\varepsilon }=\frac{1}{\alpha _{\varepsilon }\beta
_{\varepsilon }}$ and let%
\[
\tilde{u}_{\varepsilon }\left( x,t\right) =\left( \sqrt{\alpha _{\varepsilon
}\beta _{\varepsilon }}\mu _{\varepsilon }\left( t\right) \right) ^{\frac{n}{%
2}}u\left( \sqrt{\alpha _{\varepsilon }\beta _{\varepsilon }}x\mu
_{\varepsilon }\left( t\right) ,\beta _{\varepsilon }t\mu _{\varepsilon
}\left( t\right) \right) e^{\eta } 
\]%
be the function associated to $u_{\varepsilon }$ in Lemma 4.1, where $%
a+ib=\varepsilon +i$ and $\alpha ,$ $\beta $ are replaced respectively by $%
\alpha _{\varepsilon }$, $\beta _{\varepsilon }$ and

\[
\mu _{\varepsilon }\left( t\right) =\frac{1}{\alpha _{\varepsilon }\left(
1-t\right) +\beta _{\varepsilon }t}. 
\]%
Because $\alpha <\beta $, $\tilde{u}_{\varepsilon }\in L^{\infty }\left(
0,1;X\right) \cap L^{2}\left( 0,1;Y^{1}\right) $ and satisfies the equation%
\[
\partial _{t}\tilde{u}_{\varepsilon }=\left( \varepsilon +i\right) \left(
\Delta \tilde{u}_{\varepsilon }+A\tilde{u}_{\varepsilon }+\tilde{V}%
_{1}^{\varepsilon }\left( x,t\right) \tilde{u}_{\varepsilon }+\tilde{F}%
_{\varepsilon }\left( x,t\right) \right) \text{ in }R^{n}\times \left[ 0,1%
\right] , 
\]%
where 
\[
\tilde{V}_{1}^{\varepsilon }\left( x,t\right) =\alpha _{\varepsilon }\beta
_{\varepsilon }\mu _{\varepsilon }\left( t\right) V_{1}\left( \sqrt{\alpha
_{\varepsilon }\beta _{\varepsilon }}\mu _{\varepsilon }\left( t\right)
x\right) ,\text{ }\sup\limits_{t\in \left[ 0,1\right] }\left\Vert
e^{\left\vert x\right\vert ^{2}\mu ^{2}\left( t\right) }\tilde{V}%
_{1}^{\varepsilon }\left( ,.t\right) \right\Vert _{B}\leq \frac{\beta }{%
\alpha }M_{1}, 
\]%
\begin{equation}
\tilde{F}_{\varepsilon }\left( x,t\right) =\left[ \sqrt{\alpha _{\varepsilon
}\beta _{\varepsilon }}\mu _{\varepsilon }\left( t\right) \right] ^{\frac{n}{%
2}+2}V_{1}\left( \sqrt{\alpha _{\varepsilon }\beta _{\varepsilon }}\mu
_{\varepsilon }\left( t\right) x,\beta _{\varepsilon }t\mu _{\varepsilon
}\left( t\right) \right) ,  \tag{5.9}
\end{equation}

\begin{equation}
\left\Vert e^{\gamma _{\varepsilon }\left\vert x\right\vert ^{2}}\tilde{F}%
_{\varepsilon }\left( .,t\right) \right\Vert _{X}\leq \frac{\beta }{\alpha }%
\left\Vert e^{\mu _{\varepsilon }^{2}\left\vert x\right\vert
^{2}}F_{\varepsilon }\left( .,t\right) \right\Vert _{X}\text{, }\left\Vert 
\tilde{F}_{\varepsilon }\left( .,t\right) \right\Vert _{X}\leq \frac{\beta }{%
\alpha }\left\Vert F_{\varepsilon }\left( .,t\right) \right\Vert _{X} 
\tag{5.10}
\end{equation}%
and 
\begin{equation}
\left\Vert e^{\gamma _{\varepsilon }\left\vert x\right\vert ^{2}}\tilde{u}%
_{\varepsilon }\left( .,t\right) \right\Vert _{X}=\left\Vert e^{\left[ \mu
_{\varepsilon }^{2}+\frac{\left( \alpha _{\varepsilon }-\beta _{\varepsilon
}\right) a}{4\left( a^{2}+b^{2}\right) }\gamma _{\varepsilon }\right]
\left\vert y\right\vert ^{2}}u_{\varepsilon }\left( .,s\right) \right\Vert
_{X},\text{ }\left\Vert \tilde{u}_{\varepsilon }\left( .,t\right)
\right\Vert _{X}\leq \left\Vert u\left( .,s\right) \right\Vert _{X}, 
\tag{5.11}
\end{equation}%
when $s=\mu _{\varepsilon }\left( t\right) .$ The above identity when $t$ is
zero or one and (5.6) shows that 
\begin{equation}
\left\Vert e^{\gamma _{\varepsilon }\left\vert x\right\vert ^{2}}\tilde{u}%
_{\varepsilon }\left( .,0\right) \right\Vert _{X}\leq \left\Vert e^{\frac{%
\left\vert x\right\vert ^{2}}{\beta ^{2}}}u\left( .,0\right) \right\Vert
_{X},\text{ }\left\Vert e^{\gamma _{\varepsilon }\left\vert x\right\vert
^{2}}\tilde{u}_{\varepsilon }\left( .,1\right) \right\Vert _{X}\leq 
\tag{5.12}
\end{equation}%
\[
e^{\varepsilon \left\Vert V\right\Vert _{B}}\left\Vert V_{2}\left(
.,t\right) \right\Vert _{B}\left\Vert e^{\frac{\left\vert x\right\vert ^{2}}{%
\beta ^{2}}}u\left( .,1\right) \right\Vert _{X}. 
\]%
On the other hand, 
\begin{equation}
N_{1}^{-1}\left\Vert u\left( .,0\right) \right\Vert _{X}\leq \left\Vert
u\left( .,t\right) \right\Vert _{X}\leq N_{1}\left\Vert u\left( .,0\right)
\right\Vert _{X},\text{ }t\in \left[ 0,1\right] ,  \tag{5.13}
\end{equation}%
where 
\[
N_{1}=e^{\sup\limits_{t\in \left[ 0,1\right] }\left\Vert \func{Re}%
V_{2}\left( .,t\right) \right\Vert _{B}}. 
\]%
The energy method imply that 
\begin{equation}
\partial _{t}\left\Vert \tilde{u}_{\varepsilon }\left( .,t\right)
\right\Vert _{X}^{2}\leq 2\varepsilon \left\Vert \tilde{V}_{1}^{\varepsilon
}\left( x,t\right) \right\Vert _{B}\left\Vert \tilde{u}_{\varepsilon }\left(
.,t\right) \right\Vert _{X}^{2}+2\left\Vert \tilde{F}_{\varepsilon }\left(
x,t\right) \right\Vert _{X}\left\Vert \tilde{u}_{\varepsilon }\left(
.,t\right) \right\Vert _{X}.  \tag{5.14}
\end{equation}

Let $0=t_{0}<t_{1}<....t_{m}=1$ be a uniformly distributed partition of $%
[0,1]$, where $m$ will be chosen later.The inequalities $(5.14),(5.9$), the
inequality in $(5.11)$, the second inequality in $(5.10),(5.8)$ and $(5.13)$
imply that there is $N_{2}$, which depends on $\frac{\beta }{\alpha }$, $%
\left\Vert V_{1}\right\Vert _{B}$ and $\sup\limits_{t\in \left[ 0,1\right]
}\left\Vert \func{Re}V_{2}\left( .,t\right) \right\Vert _{B}$ such that%
\begin{equation}
\left\Vert \tilde{u}_{\varepsilon }\left( .,t_{i}\right) \right\Vert
_{X}\leq e^{\frac{\varepsilon \beta }{\alpha }\left\Vert V_{1}\right\Vert
_{B}}\left\Vert \tilde{u}_{\varepsilon }\left( .,t\right) \right\Vert
_{X}+N_{2}\sqrt{t_{i}-t_{i-1}}\left\Vert u\left( .,0\right) \right\Vert _{X}
\tag{5.15}
\end{equation}%
for $t\in \left[ t_{i-1},t_{i}\right] $ and $i=1,2,...m.$ Choose now $m$ so
that%
\begin{equation}
N_{2}\max\limits_{i}\sqrt{t_{i}-t_{i-1}}\leq \frac{1}{4N_{1}}.  \tag{5.16}
\end{equation}%
Because, $\lim\limits_{\varepsilon \rightarrow 0}\left\Vert \tilde{u}%
_{\varepsilon }\left( .,t\right) \right\Vert _{X}=\left\Vert u\left(
.,s\right) \right\Vert _{X}$ when $s=\beta t\mu \left( t\right) $ and $%
(5.13) $, there is $\varepsilon _{0}$ such that 
\begin{equation}
\left\Vert \tilde{u}_{\varepsilon }\left( .,t_{i}\right) \right\Vert
_{X}\geq \frac{1}{4N_{1}}\left\Vert u\left( .,0\right) \right\Vert _{X},%
\text{ when }0<\varepsilon \leq \varepsilon _{0},\text{ }i=1,2,...m 
\tag{5.17}
\end{equation}%
and now, $(5.15)-(5.17)$ show that 
\begin{equation}
\left\Vert \tilde{u}_{\varepsilon }\left( .,t\right) \right\Vert _{X}\geq 
\frac{1}{4N_{1}}\left\Vert u\left( .,0\right) \right\Vert _{X},\text{ when }%
0<\varepsilon \leq \varepsilon _{0},\text{ }t\in \left[ 0,1\right] . 
\tag{5.18}
\end{equation}%
It is now simple to verify that $(5.18)$, the first inequality in $(5.10),$ $%
(5.7)$ and $(5.13)$ imply that%
\begin{equation}
\sup\limits_{t\in \left[ 0,1\right] }\frac{\left\Vert e^{\gamma
_{\varepsilon }\left\vert x\right\vert ^{2}}\tilde{F}_{\varepsilon }\left(
.,t\right) \right\Vert _{X}}{\left\Vert \tilde{u}_{\varepsilon }\left(
.,t\right) \right\Vert _{X}}\leq \frac{4\beta }{\alpha }M_{2}\left(
\varepsilon \right) ,\text{ when }0<\varepsilon \leq \varepsilon _{0} 
\tag{5.19}
\end{equation}%
where 
\[
M_{2}\left( \varepsilon \right) =e^{\sup\limits_{t\in \left[ 0,1\right]
}\left\Vert \func{Re}V_{2}\left( .,t\right) \right\Vert _{B}+\varepsilon
\left\Vert V_{1}\right\Vert _{B}}\sup\limits_{t\in \left[ 0,1\right]
}\left\Vert e^{\left\vert x\right\vert ^{2}\mu ^{2}\left( t\right)
}V_{2}\left( .,t\right) \right\Vert _{B}. 
\]%
By using Lemma $3$.$3,$ $(5.12),$ $(5.9)$ and $(5.19)$ to show that $%
\left\Vert e^{\gamma _{\varepsilon }\left\vert x\right\vert ^{2}}\tilde{u}%
_{\varepsilon }\left( .,t\right) \right\Vert _{X}$ is logarithmically convex
in $[0,1]$ and that%
\begin{equation}
\left\Vert e^{\gamma \left\vert x\right\vert ^{2}}\tilde{u}_{\varepsilon
}\left( .,t\right) \right\Vert _{X}\leq e^{NM\left( a,b\right) }\left\Vert
e^{\gamma \left\vert x\right\vert ^{2}}\tilde{u}_{\varepsilon }\left(
0\right) \right\Vert _{X}^{1-t}\left\Vert e^{\gamma \left\vert x\right\vert
^{2}}\tilde{u}_{\varepsilon }\left( 1\right) \right\Vert _{X}^{t}, 
\tag{5.20}
\end{equation}%
when $0<\varepsilon \leq \varepsilon _{0},$ $t\in \left[ 0,1\right] $ and $%
N=N\left( \alpha ,\beta \right) .$ Then, Lemma 3.4 gives that%
\[
\left\Vert \eta \nabla \tilde{u}_{\varepsilon }\right\Vert _{Z}+\left\Vert
\eta \left\vert x\right\vert \tilde{u}_{\varepsilon }\right\Vert _{Z}\leq 
\]

\[
N\left( 1+M_{1}\right) \left[ \sup\limits_{t\in \left[ 0,1\right]
}\left\Vert e^{\gamma \left\vert x\right\vert ^{2}}\tilde{u}_{\varepsilon
}\left( .,t\right) \right\Vert _{X}+\sup\limits_{t\in \left[ 0,1\right]
}\left\Vert e^{\gamma \left\vert x\right\vert ^{2}}\tilde{F}_{\varepsilon
}\left( .,t\right) \right\Vert _{Z}\right] \leq 
\]%
\[
Ne^{N\left( M_{0}+M_{1}+M_{2}\left( \varepsilon \right)
+M_{1}^{2}+M_{2}^{2}\left( \varepsilon \right) \right) }\left[ \left\Vert e^{%
\frac{\left\vert x\right\vert ^{2}}{\beta ^{2}}}u\left( .,0\right)
\right\Vert _{X}+\left\Vert e^{\frac{\left\vert x\right\vert ^{2}}{\alpha
^{2}}}u\left( .,1\right) \right\Vert _{X}\right] , 
\]%
when $0<\varepsilon \leq \varepsilon _{0}$, the logarithmic convexity and
regularity of $u$ follow from the limit of the identity in $(5.11)$, the
final limit relation between the variables $s$ and $t$, $s=$ $\beta t\mu
\left( t\right) $ and letting $\varepsilon $ tend to zero in $(5.20)$ and
the above inequality.

By reasoning as in $\left[ \text{4, Lemma 6}\right] $ we obtain:

\textbf{Lemma 5.1.} Assume $A$ is a symmetric operator in $H$ and $V\left(
x,t\right) $ is a potential operator function in $H$ such that 
\[
\left\Vert V\right\Vert _{B}\leq \varepsilon _{0}\text{ for a }\varepsilon
_{0}>0. 
\]

Let $u\in C\left( \left[ 0,1\right] ;X\right) $ be a solution of the
equation 
\[
\partial _{t}u=i\left[ \Delta u+Au+V\left( x,t\right) u+F\left( x,t\right) %
\right] ,\text{ }x\in R^{n},\text{ }t\in \left[ 0,1\right] . 
\]

\bigskip Then, 
\[
\sup\limits_{t\in \left[ 0,1\right] }\left\Vert e^{\lambda .x}u\left(
.,t\right) \right\Vert _{X}\leq N\left[ \left\Vert e^{\lambda .x}u\left(
.,0\right) \right\Vert _{X}+\left\Vert e^{\lambda .x}u\left( .,1\right)
\right\Vert _{X}+\left\Vert e^{\lambda .x}F\left( .,t\right) \right\Vert
_{L^{1}\left( 0,1;X\right) }\right] , 
\]%
where $\lambda \in R^{n}$ and $N>0$ is constant.

\textbf{Theorem 5.1.} Assume $A$ is a symmetric operator in $H$ and $V\left(
x,t\right) $ is a potential operator function in $H$ such that%
\[
V\in B\text{ and }\lim\limits_{R\rightarrow \infty }\left\Vert V\right\Vert
_{B\left( R\right) }=0. 
\]

Suppose $\alpha ,$ $\beta $ are positive numbers and 
\[
\left\Vert e^{\frac{\left\vert x\right\vert ^{2}}{\beta ^{2}}}u\left(
.,0\right) \right\Vert _{X}<\infty \text{, }\left\Vert e^{\frac{\left\vert
x\right\vert ^{2}}{\alpha ^{2}}}u\left( .,1\right) \right\Vert _{X}<\infty . 
\]

Let $u\in C\left( \left[ 0,1\right] ;X\right) $ be a solution of the
equation 
\[
\partial _{t}u=i\left[ \Delta u+Au+V\left( x,t\right) u\right] ,\text{ }x\in
R^{n},\text{ }t\in \left[ 0,1\right] . 
\]

\bigskip Then, there is a $N=N(\alpha ,\beta )$ such that 
\[
\sup\limits_{t\in \left[ 0,1\right] }\left\Vert e^{\left\vert x\right\vert
^{2}\mu ^{2}\left( t\right) }u\left( .,t\right) \right\Vert _{X}+\left\Vert 
\sqrt{t\left( 1-t\right) }e^{\left\vert x\right\vert ^{2}\mu ^{2}\left(
t\right) }\nabla u\right\Vert _{L^{2}\left( R^{n}\times \left[ 0,1\right]
;H\right) }\leq 
\]%
\[
Ne^{\sup\limits_{t\in \left[ 0,1\right] }\left\Vert V\right\Vert _{B}}\left[
\left\Vert e^{\frac{\left\vert x\right\vert ^{2}}{\beta ^{2}}}u\left(
.,0\right) \right\Vert _{X}+\left\Vert e^{\frac{\left\vert x\right\vert ^{2}%
}{\alpha ^{2}}}u\left( .,1\right) \right\Vert _{X}+\sup\limits_{t\in \left[
0,1\right] }\left\Vert u\left( .,t\right) \right\Vert _{X}\right] . 
\]

\textbf{Proof. }Assume that $u(y,s)$ verifies the equation%
\[
\partial _{s}u=i\left[ \Delta u+Au+V\left( y,s\right) u+F\left( y,s\right) %
\right] ,\text{ }y\in R^{n},\text{ }s\in \left[ 0,1\right] . 
\]%
Set $\gamma =\left( \alpha \beta \right) ^{-1}$ and let 
\begin{equation}
\tilde{u}\left( x,t\right) =\left( \sqrt{\alpha \beta }\mu \left( t\right)
\right) ^{\frac{n}{2}}u\left( \sqrt{\alpha \beta }x\mu \left( t\right)
,\beta t\mu \left( t\right) \right) e^{\eta }.  \tag{5.21}
\end{equation}%
The function $\left( 5.21\right) $ is a solution of 
\[
\partial _{t}u=i\left[ \Delta u+Au+V\left( x,t\right) u\right] ,\text{ }x\in
R^{n},\text{ }t\in \left[ 0,1\right] 
\]%
with 
\[
\tilde{V}\left( x,t\right) =\alpha \beta \mu ^{2}\left( t\right) V\left( 
\sqrt{\alpha \beta }x\mu \left( t\right) ,\beta t\mu \left( t\right) \right)
, 
\]%
\[
\sup\limits_{t\in \left[ 0,1\right] }\left\Vert \tilde{V}\left( .,t\right)
\right\Vert _{B}\leq \max \left( \frac{\alpha }{\beta },\frac{\beta }{\alpha 
}\right) \sup\limits_{t\in \left[ 0,1\right] }\left\Vert V\left( .,t\right)
\right\Vert _{B},\text{ }\lim\limits_{R\rightarrow \infty }\left\Vert \tilde{%
V}\left( .,t\right) \right\Vert _{B\left( R\right) }=0 
\]%
and%
\begin{equation}
\left\Vert e^{\gamma \left\vert x\right\vert ^{2}}\tilde{u}\left( .,t\right)
\right\Vert _{X}=\left\Vert e^{\mu ^{2}\left( t\right) \left\vert
x\right\vert ^{2}}u\left( .,s\right) \right\Vert _{X},  \tag{5.22 }
\end{equation}%
\[
\left\Vert \tilde{u}\left( .,t\right) \right\Vert _{X}=\left\Vert u\left(
.,s\right) \right\Vert _{X}\text{ when }s=\beta t\mu \left( t\right) . 
\]%
Choose $R>0$ such that $\left\Vert \tilde{V}\left( .,t\right) \right\Vert
_{B\left( R\right) }\leq \varepsilon _{0}$ we get%
\[
\partial _{t}\tilde{u}=i\left[ \Delta \tilde{u}+A\tilde{u}+\tilde{V}%
_{R}\left( x,t\right) u+\tilde{F}_{R}\left( x,t\right) \right] ,\text{ }x\in
R^{n},\text{ }t\in \left[ 0,1\right] , 
\]%
with 
\[
\tilde{V}_{R}\left( x,t\right) =\chi _{R^{n}/O_{R}}\tilde{V}\left(
x,t\right) \text{, }\tilde{F}_{R}\left( x,t\right) =\chi _{O_{R}}\tilde{V}%
\left( x,t\right) \tilde{u}. 
\]
Then using the Lemma 5.1 we obtain 
\[
\sup\limits_{t\in \left[ 0,1\right] }\left\Vert e^{\lambda .x}\tilde{u}%
\left( .,t\right) \right\Vert _{X}\leq 
\]%
\[
N\left[ \left\Vert e^{\lambda .x}\tilde{u}\left( .,0\right) \right\Vert
_{X}+\left\Vert e^{\lambda .x}\tilde{u}\left( .,1\right) \right\Vert
_{X}+e^{\left\vert \lambda \right\vert R}\left\Vert \tilde{V}\left(
.,t\right) \right\Vert _{B}\sup\limits_{t\in \left[ 0,1\right] }\left\Vert
u\left( .,t\right) \right\Vert _{X}\right] . 
\]%
Replace $\lambda $ by $\lambda \sqrt{\gamma }$ in the above inequality,
square both sides, multiply all by $e^{-\frac{\left\vert \lambda
^{2}\right\vert }{2}}$ and integrate both sides with respect to $\lambda $
in $R^{n}$. This and the identity%
\[
\dint\limits_{R^{n}}e^{2\sqrt{\gamma }\lambda .x-\frac{\left\vert \lambda
\right\vert ^{2}}{2}}d\lambda =\left( 2\pi \right) ^{\frac{n}{2}}e^{2\gamma
\left\vert x\right\vert ^{2}} 
\]%
imply the inequality 
\begin{equation}
\sup\limits_{t\in \left[ 0,1\right] }\left\Vert \tilde{u}\left( .,t\right)
\right\Vert _{X}\leq  \tag{5.23}
\end{equation}%
\[
N\left[ \left\Vert e^{2\gamma \left\vert x\right\vert ^{2}}\tilde{u}\left(
.,0\right) \right\Vert _{X}+\left\Vert e^{2\gamma \left\vert x\right\vert
^{2}}\tilde{u}\left( .,1\right) \right\Vert _{X}+\left\Vert e^{2\gamma R^{2}}%
\tilde{V}\left( .,t\right) \right\Vert _{B}\sup\limits_{t\in \left[ 0,1%
\right] }\left\Vert \tilde{u}\left( .,t\right) \right\Vert _{X}\right] . 
\]%
This inequality and $(5.22)$ imply that%
\[
\sup\limits_{t\in \left[ 0,1\right] }\left\Vert \tilde{u}\left( .,t\right)
\right\Vert _{X}\leq 
\]%
\[
N\left[ \left\Vert e^{\frac{\left\vert x\right\vert ^{2}}{\beta ^{2}}}\tilde{%
u}\left( .,0\right) \right\Vert _{X}+\left\Vert e^{\frac{\left\vert
x\right\vert ^{2}}{\beta ^{2}}}\tilde{u}\left( .,1\right) \right\Vert
_{X}+\sup\limits_{t\in \left[ 0,1\right] }\left\Vert V\left( .,t\right)
\right\Vert _{B}\sup\limits_{t\in \left[ 0,1\right] }\left\Vert u\left(
.,t\right) \right\Vert _{X}\right] 
\]%
for some new constant $N$.

To prove the regularity of $u$ we proceed as in $(5.2)-$ $(5.4)$. The
Duhamel formula shows that%
\begin{equation}
u_{\varepsilon }\left( x,t\right)
=e^{itW}u_{0}+i\dint\limits_{0}^{t}e^{i\left( t-s\right) W}V_{2}\left(
x,s\right) u\left( x,s\right) ds,\text{ }x\in R^{n},\text{ }t\in \left[ 0,1%
\right] .  \tag{5.24}
\end{equation}%
For $0\leq \varepsilon \leq 1$, set 
\begin{equation}
\tilde{F}_{\varepsilon }\left( x,t\right) =\frac{i}{\varepsilon +i}%
e^{\varepsilon t\left( \Delta +A\right) }\tilde{V}\left( x,t\right) \tilde{u}%
\left( x,t\right) ,  \tag{5.25}
\end{equation}%
and 
\begin{equation}
\tilde{u}_{\varepsilon }\left( x,t\right) =e^{\left( \varepsilon +i\right)
t\left( \Delta +A\right) }u_{0}+\left( \varepsilon +i\right)
\dint\limits_{0}^{t}e^{\left( \varepsilon +i\right) \left( t-s\right) \left(
\Delta +A\right) }\tilde{F}\left( x,s\right) u\left( x,s\right) ds, 
\tag{5.26}
\end{equation}%
\[
\text{ }x\in R^{n},\text{ }t\in \left[ 0,1\right] . 
\]%
The identities 
\[
e^{\left( z_{1}+z_{2}\right) \left( \Delta +A\right) }=e^{z_{1}\left( \Delta
+A\right) }.e^{z_{2}\left( \Delta +A\right) }\text{ for }\func{Re}z_{1},%
\text{ }\func{Re}z_{2}\geq 0 
\]%
and $\left( 5.24\right) -\left( 5.26\right) $ show that%
\begin{equation}
\tilde{u}_{\varepsilon }\left( x,t\right) =e^{\varepsilon t\left( \Delta
+A\right) }\tilde{u}\left( x,t\right) \text{ for }t\in \left[ 0,1\right] . 
\tag{5.27}
\end{equation}%
From Lemma 3.1 with $a+ib=\varepsilon $, $\left( 5.27\right) $ and $\left(
5.25\right) $ we get that 
\begin{equation}
\sup\limits_{t\in \left[ 0,1\right] }\left\Vert e^{\gamma _{\varepsilon
}\left\vert x\right\vert ^{2}}\tilde{u}_{\varepsilon }\left( .,t\right)
\right\Vert _{X}\leq \sup\limits_{t\in \left[ 0,1\right] }\left\Vert
e^{\gamma \left\vert x\right\vert ^{2}}\tilde{u}\left( .,t\right)
\right\Vert _{X},  \tag{4.28}
\end{equation}%
\[
\sup\limits_{t\in \left[ 0,1\right] }\left\Vert e^{\gamma _{\varepsilon
}\left\vert x\right\vert ^{2}}\tilde{F}_{\varepsilon }\left( .,t\right)
\right\Vert _{X}\leq e^{\tilde{V}_{0}}\sup\limits_{t\in \left[ 0,1\right]
}\left\Vert e^{\gamma \left\vert x\right\vert ^{2}}\tilde{F}\left(
.,t\right) \right\Vert _{X}, 
\]%
where 
\[
\gamma _{\varepsilon }=\frac{\gamma }{1+4\gamma \varepsilon },\text{ }\tilde{%
V}_{0}=\sup\limits_{t\in \left[ 0,1\right] }\left\Vert \tilde{V}\right\Vert
_{B}. 
\]%
Then, Lemma 3.4, $(5.28)$ and $(5.23)$ show that 
\[
\left\Vert e^{\gamma _{\varepsilon }\left\vert x\right\vert ^{2}}u\left(
.,t\right) \right\Vert _{L^{2}\left( R^{n}\times \left[ 0,1\right] ;H\right)
}+\left\Vert \sqrt{t\left( 1-t\right) }e^{\gamma _{\varepsilon }\left\vert
x\right\vert ^{2}}\nabla u\right\Vert _{L^{2}\left( R^{n}\times \left[ 0,1%
\right] ;H\right) }\leq 
\]%
\[
Ne^{NV_{0}}\left[ \left\Vert e^{\frac{\left\vert x\right\vert ^{2}}{\beta
^{2}}}u\left( .,0\right) \right\Vert _{X}+\left\Vert e^{\frac{\left\vert
x\right\vert ^{2}}{\alpha ^{2}}}u\left( .,1\right) \right\Vert
_{X}+\sup\limits_{t\in \left[ 0,1\right] }\left\Vert u\left( .,t\right)
\right\Vert _{X}\right] , 
\]%
where 
\[
V_{0}=\sup\limits_{t\in \left[ 0,1\right] }\left\Vert V\left( x,t\right)
\right\Vert _{B}. 
\]%
The Theorem 5.1 follows from this inequality, from $(5.21)-(5.23)$ and
letting $\varepsilon $ tend to zero.

\begin{center}
\textbf{6. A Hardy type abstract uncertainty principle. Proof of Theorem 1.}
\end{center}

\bigskip The assertion about the Carleman inequality in Lemma 6.1 below is
the following monotonicity or frequency function argument related to Lemma
3. 2. When $u\in C([0,1];X)$ is a free solution to the free abstract Schr%
\"{o}dinger equation%
\[
\partial _{t}u-i\left( \Delta u+Au\right) =0,\text{ }x\in R^{n},\text{ }t\in %
\left[ 0,1\right] , 
\]%
satisfies

\[
\left\Vert e^{\gamma \left\vert x\right\vert ^{2}}u\left( .,0\right)
\right\Vert _{X}<\infty \text{, }\left\Vert e^{\gamma \left\vert
x\right\vert ^{2}}u\left( .,1\right) \right\Vert _{X}<\infty 
\]%
and 
\[
\ f=e^{\varkappa }u,\text{ }Q\left( t\right) =\left( f\left( .,t\right)
,f\left( .,t\right) \right) _{X}, 
\]%
where

\[
\varkappa \left( x,t\right) =%
%TCIMACRO{\U{3bc} }%
%BeginExpansion
\mu
%EndExpansion
|x+Rt(1-t)|^{2}-\frac{R^{2}t(1-t)}{8\mu },\text{ }\sigma \left( \varepsilon
,t\right) =\frac{\left( 1+\varepsilon \right) t(1-t)}{16\mu }. 
\]%
Then, $\log Q\left( t\right) $ is logaritmicaly convex in $[0,1],$ when $%
0<\mu <\gamma .$

The formal application of the above argument to a $C([0,1];X)$ solution of
the equation%
\begin{equation}
\partial _{t}u-i\left[ \Delta u+Au+V\left( x,t\right) u\right] =0,\text{ }%
x\in R^{n},\text{ }t\in \left[ 0,1\right] ,  \tag{6.1}
\end{equation}%
implies a similar result, when $V$ is a bounded potential, though the
justification of the correctness of the assertions involved in the
corresponding formal application of Lemma 3.2 were formal. In fact, we can
only justify these assertions, when the potential $V$ verifies the first
condition in Theorem 1 or when we can obtain the additional regularity of
the gradient of $u$ in the strip, as in Theorem 5.1. Here, we choose to
prove Theorem 1 using the Carleman inequality in Lemma 6.1 in place of the
above convexity argument. The reason for our choice is that it is simpler to
justify the correctness of the application of the Carleman inequality to a $%
C([0,1];X)$ solution to $(6.1)$ than the corresponding monotonicity or
logarithmic convexity of the solution.

\textbf{Lemma 6.1. }Assume $A$ is a symmetric operator in $H$ and $V\left(
x,t\right) $ is a potential operator function in $H$ such that%
\[
V\in B\text{ and }\lim\limits_{R\rightarrow \infty }\left\Vert V\right\Vert
_{B\left( R\right) }=0. 
\]

The estimate

\[
R\sqrt{\frac{\varepsilon }{8\mu }}\left\Vert e^{\varkappa -\sigma }\upsilon
\right\Vert _{L^{2}\left( R^{n+1};H\right) }\leq \left\Vert e^{\varkappa
-\sigma }\left[ \partial _{t}u-i\left( \Delta u+Au\right) \right] \upsilon
\right\Vert _{L^{2}\left( R^{n+1};H\right) } 
\]%
holds, when $\varepsilon >0$, $\mu >0$, $R>0$ and $\upsilon \in
C_{0}^{\infty }\left( R^{n+1};H\right) $.

\textbf{Proof.} Let\textbf{\ }$f=e^{\varkappa -\sigma }\upsilon $. Then, 
\[
e^{\varkappa -\sigma }\left[ \partial _{t}u-i\left( \Delta u+Au\right) %
\right] \upsilon =\partial _{t}f+Sf-Kf. 
\]%
From $(3.8)-(3.10)$ with $\gamma =1$, $a+ib=i$ and $\varphi \left(
x,t\right) =\varkappa \left( x,t\right) -\sigma \left( \varepsilon ,t\right) 
$ we have 
\[
S=-4%
%TCIMACRO{\U{3bc} }%
%BeginExpansion
\mu
%EndExpansion
i(x+Rt(1-t)e_{1})\text{\textperiodcentered }\nabla -2%
%TCIMACRO{\U{3bc} }%
%BeginExpansion
\mu
%EndExpansion
ni+2%
%TCIMACRO{\U{3bc} }%
%BeginExpansion
\mu
%EndExpansion
R(1-2t)(x_{1}+Rt(1-t))-\sigma ,\text{ } 
\]%
\[
K=i\left( \triangle +A\right) +4%
%TCIMACRO{\U{3bc} }%
%BeginExpansion
\mu
%EndExpansion
^{2}i|x+Rt(1-t)e_{1}|^{2}\text{, }S_{t}+[S,K]= 
\]%
\[
-8%
%TCIMACRO{\U{3bc} }%
%BeginExpansion
\mu
%EndExpansion
\triangle +32%
%TCIMACRO{\U{3bc} }%
%BeginExpansion
\mu
%EndExpansion
^{3}|x+Rt(1-t)e_{1}|^{2}-4%
%TCIMACRO{\U{3bc} }%
%BeginExpansion
\mu
%EndExpansion
R(x_{1}+Rt(1-t))+ 
\]

\[
+2%
%TCIMACRO{\U{3bc} }%
%BeginExpansion
\mu
%EndExpansion
R^{2}(1-2t)^{2}+\frac{\left( 1+\varepsilon \right) R^{2}}{8\mu }+-4i%
%TCIMACRO{\U{3bc} }%
%BeginExpansion
\mu
%EndExpansion
R(1-2t)\partial _{x_{1}} 
\]%
and

\begin{equation}
(S_{t}f+[S,K]f,f)_{X}=32%
%TCIMACRO{\U{3bc} }%
%BeginExpansion
\mu
%EndExpansion
^{3}\dint\limits_{R^{n}}\left\vert x+Rt(1-t)e_{1}-\frac{R}{16\mu ^{2}}%
e_{1}\right\vert ^{2}\left\Vert f\right\Vert _{H}^{2}dx+  \tag{6.2}
\end{equation}%
\[
\frac{\varepsilon R^{2}}{8\mu }\dint\limits_{R^{n}}\left\Vert f\right\Vert
_{H}^{2}dx+8\mu \dint\limits_{R^{n}}\left\Vert \nabla _{x^{\prime
}}f\right\Vert _{H}^{2}dx+8\mu \dint\limits_{R^{n}}\left\Vert i\partial
_{x_{1}}f-R\left( \frac{1}{2}-t\right) f\right\Vert _{H}^{2}dx\geq 
\]%
\[
\frac{\varepsilon R^{2}}{8\mu }\dint\limits_{R^{n}}\left\Vert f\right\Vert
_{H}^{2}dx. 
\]%
Following the standard method to handle $L_{2}$-Carleman inequalities, the
symmetric and skew-symmetric parts of $\partial _{t}-S-K$, as a space-time
operator, are respectively $-S$ and $\partial _{t}-K$, and $[-S,\partial
_{t}-K]=S_{t}+[S,K]$. Thus, 
\[
\left\Vert \partial _{t}f-Sf-Kf\right\Vert _{L^{2}\left( R^{n+1};H\right)
}^{2}=\left\Vert \partial _{t}f-Kf\right\Vert _{L^{2}\left( R^{n+1};H\right)
}^{2}+\left\Vert Sf\right\Vert _{L^{2}\left( R^{n+1};H\right) }^{2}- 
\]%
\begin{equation}
2\func{Re}\dint\limits_{R^{n}}\dint\limits_{-\infty }^{\infty }\left(
Sf,\partial _{t}f-Kf\right) _{H}dxdt\geq
\dint\limits_{R^{n}}\dint\limits_{-\infty }^{\infty }\left( [-S,\partial
_{t}-K]f,f\right) _{H}dxdt=  \tag{6.3}
\end{equation}%
\[
\dint\limits_{-\infty }^{\infty }\left( S_{t}f+\left[ S,K\right] f,f\right)
_{H}dt, 
\]%
and the Lemma 6.1 follows from $(5.2)$ and $(5.3)$.

\textbf{Proof of Theorem 1. }Let u be as in Theorem $1$ and $\tilde{u}$, $%
\tilde{V}$ the corresponding functions defined in Lemma 4.1, when $a+ib=i$.
Then, $\tilde{u}\in $ $C([0,1];X)$ is a solution of the equation%
\[
\partial _{t}u-i\left[ \Delta u+Au+\tilde{V}u\right] =0,\text{ }x\in R^{n},%
\text{ }t\in \left[ 0,1\right] 
\]%
and 
\[
\left\Vert e^{\gamma \left\vert x\right\vert ^{2}}\tilde{u}\left( .,0\right)
\right\Vert _{X}<\infty \text{, }\left\Vert e^{\gamma \left\vert
x\right\vert ^{2}}\tilde{u}\left( .,1\right) \right\Vert _{X}<\infty \text{
for }\gamma =\frac{1}{\alpha \beta },\text{ }\gamma >\frac{1}{2}. 
\]%
The proofs of Theorem 3 show that in either case%
\begin{equation}
N_{\gamma }=\sup\limits_{t\in \left[ 0,1\right] }\left[ \left\Vert e^{\gamma
_{\varepsilon }\left\vert x\right\vert ^{2}}\tilde{u}\left( .,t\right)
\right\Vert _{L^{2}\left( R^{n}\times \left[ 0,1\right] ;H\right)
}+\left\Vert \sqrt{t\left( 1-t\right) }e^{\gamma _{\varepsilon }\left\vert
x\right\vert ^{2}}\nabla \tilde{u}\right\Vert _{L^{2}\left( R^{n}\times %
\left[ 0,1\right] ;H\right) }\right] <\infty .  \tag{6.4}
\end{equation}%
For given $R>0$, choose $%
%TCIMACRO{\U{3bc} }%
%BeginExpansion
\mu
%EndExpansion
$ and $\varepsilon $ such that%
\begin{equation}
\frac{\left( 1+\varepsilon \right) ^{\frac{3}{2}}}{2\left( 1-\varepsilon
\right) ^{3}}<\mu \leq \frac{\gamma }{1+\varepsilon }  \tag{6.5}
\end{equation}%
and let $\eta _{M}$ and $\theta _{R}$ be smooth functions verifying, $\theta
_{M}\left( x\right) $ $=1$, when $|x|\leq M$, $\theta _{M}\left( x\right) $ $%
=0$, when $|x|>2M$, $M\geq 2R$, $\eta _{R}\in C_{0}^{\infty }(0,1),$ $0\leq
\eta _{R}\left( t\right) \leq 1$, $\eta _{R}\left( t\right) =1$ for $t\in %
\left[ \frac{1}{R},1-\frac{1}{R}\right] $ and $\eta _{R}\left( t\right) =0$
for $t\in \left[ 0,\frac{1}{2R}\right] \cup \left[ 1-\frac{1}{2R},1\right] .$
Then, $\upsilon \left( x,t\right) =\eta _{R}\left( t\right) \theta
_{M}\left( x\right) \tilde{u}\left( x,t\right) $ is compactly supported in $%
R^{n}\times (0,1)$ and 
\begin{equation}
\partial _{t}\upsilon -i\left[ \Delta \upsilon +A\upsilon +\tilde{V}\upsilon %
\right] =\eta _{R}^{\prime }\left( t\right) \theta _{M}\left( x\right) 
\tilde{u}\left( x,t\right) -\left( 2\nabla \theta _{M}.\nabla \tilde{u}+%
\tilde{u}\Delta \theta _{M}\right) \eta _{R}.\text{ }  \tag{6.6}
\end{equation}%
The terms on the right hand side of $(6.6)$ are supported, where 
\[
%TCIMACRO{\U{3bc} }%
%BeginExpansion
\mu
%EndExpansion
|x+Rt(1-t)|^{2}\leq \gamma \left\vert x\right\vert ^{2}+\frac{\gamma }{%
\varepsilon }, 
\]%
\[
%TCIMACRO{\U{3bc} }%
%BeginExpansion
\mu
%EndExpansion
|x+Rt(1-t)e_{1}|^{2}\leq \gamma \left\vert x\right\vert ^{2}+\frac{\gamma }{%
\varepsilon }R^{2}. 
\]%
Apply now Lemma 6.1 to $\upsilon $ with the values of $%
%TCIMACRO{\U{3bc} }%
%BeginExpansion
\mu
%EndExpansion
$ and $\varepsilon $ chosen in $(6.5)$. This, the bounds for $%
%TCIMACRO{\U{3bc} }%
%BeginExpansion
\mu
%EndExpansion
|x+Rt(1-t)e_{1}|^{2}$ in each of the parts of the support of%
\[
\partial _{t}\upsilon -i\left[ \Delta \upsilon +A\upsilon +\tilde{V}\upsilon %
\right] 
\]%
and the natural bounds for $\nabla \theta _{M}$, $\triangle \theta _{M}$ and 
$\eta _{R}^{\prime }$ show that there is a constant $N_{\varepsilon }$ such
that 
\[
R\left\Vert e^{\varkappa -\sigma }\upsilon \right\Vert _{L^{\infty }\left(
R^{n}\times \left[ 0,1\right] ;H\right) }\leq 
\]%
\begin{equation}
N_{\varepsilon }\left\Vert \tilde{V}\right\Vert _{B}\left\Vert e^{\varkappa
-\sigma }\upsilon \right\Vert _{L^{2}\left( R^{n}\times \left[ 0,1\right]
;H\right) }+N_{\varepsilon }Re^{\frac{\gamma }{\varepsilon }%
}\sup\limits_{t\in \left[ 0,1\right] }\left\Vert e^{\gamma \left\vert
x\right\vert ^{2}}\tilde{u}\left( .,t\right) \right\Vert _{X}+  \tag{6.7}
\end{equation}%
\[
N_{\varepsilon }M^{-1}e^{\frac{\gamma }{\varepsilon }R^{2}}\left\Vert
e^{\gamma \left\vert x\right\vert ^{2}}\left( \left\Vert \tilde{u}%
\right\Vert _{H}+\left\Vert \nabla \tilde{u}\right\Vert _{H}\right)
\right\Vert _{L^{2}\left( R^{n}\times \sigma _{R}\right) }, 
\]%
where 
\[
\sigma _{R}=\left[ \frac{1}{2R},1-\frac{1}{2R}\right] . 
\]%
The first term on the right hand side of $(6.7)$ can be hidden in the left
hand side, when $R\geq 2N_{\varepsilon }\left\Vert \tilde{V}\right\Vert _{B}$%
, while the last tends to zero, when $M$ tends to infinity by $(6.4)$. This
and the fact that $\upsilon =\tilde{u}$ in $O_{\frac{\varepsilon \left(
1-\varepsilon ^{2}\right) ^{2}R}{4}}\times \left[ \frac{1-\varepsilon }{2},%
\frac{1+\varepsilon }{2}\right] ,$ where 
\[
\varkappa -\sigma \geq \frac{R^{2}}{16\mu }\left( 4\mu ^{2}\left(
1-\varepsilon \right) ^{6}-\left( 1+\varepsilon \right) ^{3}\right) 
\]%
and $\left( 6.5\right) $ show that%
\begin{equation}
e^{C\left( \gamma ,\varepsilon \right) }\left\Vert \tilde{u}\right\Vert
_{L^{2}\left( O_{\frac{R}{8}}\times \left[ \frac{1-\varepsilon }{2},\frac{%
1+\varepsilon }{2}\right] ;H\right) }\leq N_{\gamma ,\varepsilon }, 
\tag{6.8}
\end{equation}%
when $R\geq 2N_{\varepsilon }\left\Vert \tilde{V}\right\Vert _{B}.$ At the
same time 
\begin{equation}
N^{-1}\left\Vert \tilde{u}\left( .,0\right) \right\Vert _{X}\leq \left\Vert 
\tilde{u}\left( .,t\right) \right\Vert _{X}\leq N\left\Vert \tilde{u}\left(
.,1\right) \right\Vert _{X}  \tag{6.9}
\end{equation}%
for $0\leq t\leq 1$ and $N=e^{\sup\limits_{t\in \left[ 0,1\right]
}\left\Vert \tilde{V}\right\Vert _{B}}.$ Moreover, from $\left( 6.4\right) $
we get 
\begin{equation}
\left\Vert \tilde{u}\left( .,t\right) \right\Vert _{X}\leq N\left\Vert 
\tilde{u}\left( .,t\right) \right\Vert _{L^{2}\left( O_{\frac{R}{8}%
};H\right) }+e^{\frac{-\gamma R^{2}}{64}}N_{\gamma }\text{ when }0\leq t\leq
1  \tag{6.10}
\end{equation}%
Then, $(6.8)-(6.10)$ show that there is a constant $N_{\gamma ,\varepsilon
,V}$, which such that 
\[
e^{C\left( \gamma ,\varepsilon \right) R^{2}}\left\Vert \tilde{u}\left(
.,0\right) \right\Vert _{X}\leq N_{\gamma ,\varepsilon ,V}. 
\]%
For $R\rightarrow \infty $ we obtain $u\equiv 0.$

\textbf{Proof of Theorem 2.}

F\i rst of all we show the following Carleman inequality

\textbf{Lemma 6.2. }Assume $A$ is a symmetric operator in $H$ and $V\left(
x,t\right) $ is a potential operator function in $H$ such that%
\[
V\in B\text{ and }\lim\limits_{R\rightarrow \infty }\left\Vert V\right\Vert
_{B\left( R\right) }=0. 
\]

The estimate

\begin{equation}
R\sqrt{\frac{\varepsilon }{8\mu }}\left\Vert e^{\varkappa -\sigma +\chi
}\upsilon \right\Vert _{L^{2}\left( R^{n+1};H\right) }\leq \left\Vert
e^{\varkappa -\sigma +\chi }\left[ \partial _{t}u-\Delta u-Au\right]
\upsilon \right\Vert _{L^{2}\left( R^{n+1};H\right) }  \tag{6.11}
\end{equation}%
holds, when $\varepsilon >0$, $\mu >0$, $R>0$ and $\upsilon \in
C_{0}^{\infty }\left( R^{n+1};H\right) $, where%
\[
\chi \left( t\right) =\frac{R^{2}t(1-t)\left( 1-2t\right) }{6}. 
\]

\textbf{Proof.} Let\textbf{\ }$f=e^{\varkappa +\chi -\sigma }\upsilon $.
Then, 
\[
e^{\varkappa +\chi -\sigma }\left[ \partial _{t}u-\left( \Delta u+Au\right) %
\right] \upsilon =\partial _{t}f-Sf-Kf. 
\]%
From $(3.8)-(3.10)$ with $\gamma =1$, $a+ib=1$ and $\varphi \left(
x,t\right) =\varkappa \left( x,t\right) +\chi \left( t\right) -\sigma \left(
\varepsilon ,t\right) $ we have 
\[
S=\Delta +A+4%
%TCIMACRO{\U{3bc} }%
%BeginExpansion
\mu
%EndExpansion
^{2}i|x+Rt(1-t)e_{1}|^{2}+2%
%TCIMACRO{\U{3bc} }%
%BeginExpansion
\mu
%EndExpansion
ni+\text{ } 
\]%
\[
2%
%TCIMACRO{\U{3bc} }%
%BeginExpansion
\mu
%EndExpansion
R(1-2t)(x_{1}+Rt(1-t))-\sigma +\left( t^{2}-t+\frac{1}{6}\right) R^{2}, 
\]%
\[
K=-4%
%TCIMACRO{\U{3bc} }%
%BeginExpansion
\mu
%EndExpansion
(x+Rt(1-t)e_{1})\text{\textperiodcentered }\nabla -2%
%TCIMACRO{\U{3bc} }%
%BeginExpansion
\mu
%EndExpansion
n\text{, } 
\]%
\[
S_{t}+[S,K]=-8%
%TCIMACRO{\U{3bc} }%
%BeginExpansion
\mu
%EndExpansion
\triangle +32%
%TCIMACRO{\U{3bc} }%
%BeginExpansion
\mu
%EndExpansion
^{3}|x+Rt(1-t)e_{1}|^{2}+ 
\]

\[
4%
%TCIMACRO{\U{3bc} }%
%BeginExpansion
\mu
%EndExpansion
R(4\mu \left( 1-2t-1\right) \left( (x_{1}+Rt(1-t)\right) +\left( 2t-1\right)
R^{2}+\frac{\left( 1+\varepsilon \right) R^{2}}{8\mu } 
\]%
and

\begin{equation}
(S_{t}f+[S,K]f,f)_{X}=32%
%TCIMACRO{\U{3bc} }%
%BeginExpansion
\mu
%EndExpansion
^{3}\dint\limits_{R^{n}}\left\vert x+Rt(1-t)e_{1}+\frac{(4\mu \left(
1-2t-1\right) R}{16\mu ^{2}}e_{1}\right\vert ^{2}\left\Vert f\right\Vert
_{H}^{2}dx+  \tag{6.12}
\end{equation}%
\[
8\mu \dint\limits_{R^{n}}\left\Vert \nabla f\right\Vert _{H}^{2}dx+\frac{%
\varepsilon R^{2}}{8\mu }\dint\limits_{R^{n}}\left\Vert f\right\Vert
_{H}^{2}dx\geq \frac{\varepsilon R^{2}}{8\mu }\dint\limits_{R^{n}}\left\Vert
f\right\Vert _{H}^{2}dx. 
\]

Then from $\left( 6.12\right) $ a similar way as Lemma 6.1 we obtain the
estimate $\left( 6.11\right) .$

\bigskip \textbf{Proof of Theorem 4. }Assume that $u$ verifies the
conditions in Theorem 4 and let $\tilde{u}$ be the Appel transformation of $%
u $ defined in Lemma 4.1 with $a+ib=1$, $\alpha =1$ and $\beta =1+\frac{2}{%
\beta }$. $\tilde{u}\in L^{\infty }\left( 0,1;X\right) \cap L^{2}\left(
0,1;Y^{1}\right) $ is a solution of the equation 
\[
\partial _{t}u=\Delta u+Au+\tilde{V}u,\text{ }x\in R^{n},\text{ }t\in \left[
0,1\right] 
\]%
with $\tilde{V}$ a bounded potential in $R^{n}\times \lbrack 0,1]$ and $%
\gamma =\frac{1}{2\delta }.$ Then, we have 
\[
\left\Vert e^{\gamma \left\vert x\right\vert ^{2}}\tilde{u}\left( .,0\right)
\right\Vert _{X}=\left\Vert \tilde{u}\left( .,0\right) \right\Vert _{X}\text{%
, }\left\Vert e^{\gamma \left\vert x\right\vert ^{2}}\tilde{u}\left(
.,1\right) \right\Vert _{X}=\left\Vert \tilde{u}\left( .,1\right)
\right\Vert _{X}. 
\]%
From Lemma 3.3 and Lemma 3. 4 with $a+ib=1$, we have

\[
\sup\limits_{t\in \left[ 0,1\right] }\left\Vert e^{\gamma \left\vert
x\right\vert ^{2}}\tilde{u}\left( .,t\right) \right\Vert _{X}+\left\Vert 
\sqrt{t\left( 1-t\right) }e^{\gamma \left\vert x\right\vert ^{2}}\nabla 
\tilde{u}\right\Vert _{L^{2}\left( R^{n}\times \left[ 0,1\right] ;H\right)
}\leq 
\]%
\[
e^{\left( M_{1}+M_{1}^{2}\right) }\left[ \left\Vert e^{\gamma \left\vert
x\right\vert ^{2}}\tilde{u}\left( .,0\right) \right\Vert _{X}+\left\Vert
e^{\gamma \left\vert x\right\vert ^{2}}\tilde{u}\left( .,1\right)
\right\Vert _{X}\right] , 
\]%
where 
\[
M_{1}=\left\Vert \tilde{V}\right\Vert _{B}. 
\]

The proof is finished by setting $\upsilon (x,t)=\theta _{M}(x)\eta _{R}(t)%
\tilde{u}(x,t)$, by using Carleman inequality $\left( 6.11\right) $ and in\
similar argument that we used to prove Theorem 1.

\begin{center}
\textbf{7. Unique continuation properties for the system of Schr\"{o}dinger
equation }
\end{center}

Consider the Cauchy problem for the system of Schr\"{o}dinger equation%
\begin{equation}
\frac{\partial u_{m}}{\partial t}=i\left[ \Delta
u_{m}+\sum\limits_{j=1}^{N}a_{mj}u_{j}+\sum\limits_{j=1}^{N}b_{mj}u_{j}%
\right] ,\text{ }x\in R^{n},\text{ }t\in \left( 0,T\right) ,  \tag{7.1}
\end{equation}%
where $u=\left( u_{1},u_{2},...,u_{N}\right) ,$ $u_{j}=u_{j}\left(
x,t\right) ,$ $a_{mj}$ are complex numbers and $b_{mj}=b_{mj}\left(
x,t\right) $ are complex valued functions$.$ Let $l_{2}=l_{2}\left( N\right) 
$ and $l_{2}^{s}=l_{2}^{s}\left( N\right) $ (see $\left[ \text{23, \S\ 1.18}%
\right] $). Let $A$ be the operator in $l_{2}\left( N\right) $ defined by%
\[
\text{ }D\left( A\right) =\left\{ u=\left\{ u_{j}\right\} ,\text{ }%
\left\Vert u\right\Vert _{l_{2}^{s}\left( N\right) }=\left(
\sum\limits_{j=1}^{N}2^{sj}\left\vert u_{j}\right\vert ^{2}\right) ^{\frac{1%
}{2}}<\infty \right\} , 
\]

\[
A=\left[ a_{mj}\right] \text{, }a_{mj}=g_{m}2^{sj},\text{ }m,j=1,2,...,N,%
\text{ }N\in \mathbb{N} 
\]

\bigskip and 
\[
\text{ }D\left( V\left( x,t\right) \right) =\left\{ u=\left\{ u_{j}\right\}
\in l_{2}^{s}\right\} , 
\]

\[
V\left( x,t\right) =\left[ b_{mj}\left( x,t\right) \right] \text{, }%
b_{mj}\left( x,t\right) =g_{m}\left( x,t\right) 2^{sj},\text{ }%
m,j=1,2,...,N, 
\]

Let 
\[
X_{2}=L^{2}\left( R^{n};l_{2}\right) ,Y^{s,2}=H^{s,2}\left(
R^{n};l_{2}\right) . 
\]

From Theorem 1 we obtain the following result

\textbf{Theorem 7.1. }Assume $a_{mj}=a_{jm}$ and 
\[
\sup\limits_{t\in \left[ 0,1\right] }\left\Vert e^{\left\vert x\right\vert
^{2}\mu ^{2}\left( t\right) }g_{m}\left( .,t\right) \right\Vert _{L^{\infty
}\left( R^{n}\right) }<\infty . 
\]

Let $\alpha $, $\beta >0$ and $\alpha \beta <2$. Assume $u\in C\left( \left[
0,1\right] ;l_{2}\right) $ be a solution of the equation $\left( 7.1\right) $
and%
\[
\left\Vert e^{\frac{\left\vert x\right\vert ^{2}}{\beta ^{2}}}u\left(
.,0\right) \right\Vert _{X2}<\infty ,\left\Vert e^{\frac{\left\vert
x\right\vert ^{2}}{\alpha ^{2}}}u\left( .,T\right) \right\Vert
_{X_{2}}<\infty . 
\]

Then $u\left( x,t\right) \equiv 0.$

\ \textbf{Proof.} It is easy to see that $A$ is a symmetric operator in $%
l_{2}$ and other conditions of Theorem 1 are satisfied. Hence, from Teorem1
we obtain the conculision.

\begin{center}
\textbf{8. Unique continuation properties for anisotropic Schr\"{o}dinger
equation }

\ \ \ \ \ \ \ \ \ \ \ \ \ \ \ \ \ \ \ \ \ \ \ \ \ \ \ \ \ \ \ \ \ \ \ \ \ \
\ \ \ \ \ \ \ \ \ \ 
\end{center}

The regularity property of BVP for elliptic equations\ were studied e.g. in $%
\left[ \text{1, 2}\right] $. Let $\Omega =R^{n}\times G$, $G\subset R^{d},$ $%
d\geq 2$ is a bounded domain with $\left( d-1\right) $-dimensional boundary $%
\partial G$. Let us consider the following problem

\begin{equation}
\partial _{t}u=i\left[ \Delta _{x}u+\sum\limits_{\left\vert \alpha
\right\vert \leq 2m}a_{\alpha }\left( y\right) D_{y}^{\alpha }u\left(
x,y,t\right) +\dint\limits_{G}K\left( x,y,\tau ,t\right) u\left(
x,y,,t\right) d\tau \right] ,  \tag{8.1}
\end{equation}%
\[
\text{ }x\in R^{n},\text{ }y\in \Omega ,\text{ }t\in \left[ 0,T\right] ,%
\text{ } 
\]

\begin{equation}
B_{j}u=\sum\limits_{\left\vert \beta \right\vert \leq m_{j}}\ b_{j\beta
}\left( y\right) D_{y}^{\beta }u\left( x,y,t\right) =0\text{, }x\in R^{n},%
\text{ }y\in \partial G,\text{ }j=1,2,...,m,  \tag{8.2}
\end{equation}%
where $a_{\alpha },$ $a_{i\beta },$ $b_{j\beta }$ are the complex valued
functions, $\alpha =\left( \alpha _{1},\alpha _{2},...,\alpha _{n}\right) $, 
$\beta =\left( \beta _{1},\beta _{2},...,\beta _{n}\right) ,$ $\mu _{i}<2m,$ 
$K=K\left( x,y,\tau ,t\right) $ is a complex valued bounded function in $%
\Omega \times G\times \left[ 0,T\right] $ and 
\[
D_{x}^{k}=\frac{\partial ^{k}}{\partial x^{k}},\text{ }D_{j}=-i\frac{%
\partial }{\partial y_{j}},\text{ }D_{y}=\left( D_{1,}...,D_{n}\right) ,%
\text{ }y=\left( y_{1},...,y_{n}\right) ,
\]

$\bigskip $Let%
\[
\xi ^{\prime }=\left( \xi _{1},\xi _{2},...,\xi _{n-1}\right) \in R^{n-1},%
\text{ }\alpha ^{\prime }=\left( \alpha _{1},\alpha _{2},...,\alpha
_{n-1}\right) \in Z^{n},\text{ } 
\]%
\[
\text{ }A\left( y_{0},\xi ^{\prime },D_{y}\right) =\sum\limits_{\left\vert
\alpha ^{\prime }\right\vert +j\leq 2m}a_{\alpha ^{\prime }}\left(
y_{0}\right) \xi _{1}^{\alpha _{1}}\xi _{2}^{\alpha _{2}}...\xi
_{n-1}^{\alpha _{n-1}}D_{y}^{j}\text{ for }y_{0}\in \bar{G} 
\]%
\[
B_{j}\left( y_{0},\xi ^{\prime },D_{y}\right) =\sum\limits_{\left\vert \beta
^{\prime }\right\vert +j\leq m_{j}}b_{j\beta ^{\prime }}\left( y_{0}\right)
\xi _{1}^{\beta _{1}}\xi _{2}^{\beta _{2}}...\xi _{n-1}^{\beta
_{n-1}}D_{y}^{j}fory_{0}\in \partial G. 
\]

\textbf{Theorem 8.1}. Let the following conditions be satisfied:

\bigskip (1) $G\in C^{2}$, $a_{\alpha }\in C\left( \bar{G}\right) $ for each 
$\left\vert \alpha \right\vert =2m$ and $a_{\alpha }\in L_{\infty }\left(
G\right) $ for each $\left\vert \alpha \right\vert <2m$;

(2) $b_{j\beta }\in C^{2m-m_{j}}\left( \partial G\right) $ for each $j$, $%
\beta $ and $\ m_{j}<2m$, $\sum\limits_{j=1}^{m}b_{j\beta }\left( y^{\prime
}\right) \sigma _{j}\neq 0,$ for $\left\vert \beta \right\vert =m_{j},$ $%
y^{^{\shortmid }}\in \partial G,$ where $\sigma =\left( \sigma _{1},\sigma
_{2},...,\sigma _{n}\right) \in R^{n}$ is a normal to $\partial G$ $;$

(3) for $y\in \bar{G}$, $\xi \in R^{n}$, $\lambda \in S\left( \varphi
_{0}\right) $ for $0\leq \varphi _{0}<\pi $, $\left\vert \xi \right\vert
+\left\vert \lambda \right\vert \neq 0$ let $\lambda +$ $\sum\limits_{\left%
\vert \alpha \right\vert =2m}a_{\alpha }\left( y\right) \xi ^{\alpha }\neq 0$%
;

(4) for each $y_{0}\in \partial G$ local BVP in local coordinates
corresponding to $y_{0}$:%
\[
\lambda +A\left( y_{0},\xi ^{\prime },D_{y}\right) \vartheta \left( y\right)
=0, 
\]

\[
B_{j}\left( y_{0},\xi ^{\prime },D_{y}\right) \vartheta \left( 0\right)
=h_{j}\text{, }j=1,2,...,m 
\]%
has a unique solution $\vartheta \in C_{0}\left( \mathbb{R}_{+}\right) $ for
all $h=\left( h_{1},h_{2},...,h_{n}\right) \in \mathbb{C}^{n}$ and for $\xi
^{\prime }\in R^{n-1};$

(5) 
\[
\sup\limits_{t\in \left[ 0,1\right] }\left\Vert e^{\left\vert x\right\vert
^{2}\mu ^{2}\left( t\right) }K\left( .,y,t\right) \right\Vert _{L^{\infty
}\left( \Omega \times G\right) }<\infty \text{ for }y\in \bar{G}; 
\]

(6) Let $\alpha $, $\beta >0$, $\alpha \beta <2$. Assume $u\in C\left( \left[
0,1\right] ;l_{2}\right) $ be a solution of the equation $\left( 8.1\right)
-\left( 8.2\right) $ and%
\[
\left\Vert e^{\frac{\left\vert x\right\vert ^{2}}{\beta ^{2}}}u\left(
.,0\right) \right\Vert _{L^{2}\left( \Omega \right) }<\infty ,\left\Vert e^{%
\frac{\left\vert x\right\vert ^{2}}{\alpha ^{2}}}u\left( .,T\right)
\right\Vert _{L^{2}\left( \Omega \right) }<\infty . 
\]

Then $u\left( x,y,t\right) \equiv 0.$

\ \textbf{Proof. }Let us consider operators $A$ and $V\left( x,t\right) $ in 
$H=L^{2}\left( G\right) $ that are defined by the equalities 
\[
D\left( A\right) =\left\{ u\in W^{2m,2}\left( G\right) \text{, }B_{j}u=0,%
\text{ }j=1,2,...,m\text{ }\right\} ,\ Au=\sum\limits_{\left\vert \alpha
\right\vert \leq 2m}a_{\alpha }\left( y\right) D_{y}^{\alpha }u\left(
y\right) , 
\]%
\[
V\left( x,t\right) u=\dint\limits_{G}K\left( x,y,\tau ,t\right) u\left(
x,y,\tau ,t\right) d\tau . 
\]

Then the problem $\left( 8.1\right) -\left( 8.2\right) $ can be rewritten as
the problem $\left( 1.1\right) $, where $u\left( x\right) =u\left(
x,.\right) ,$ $f\left( x\right) =f\left( x,.\right) $,\ $x\in \sigma $ are
the functions with values in\ $H=L^{2}\left( G\right) $. By virtue of $\left[
\text{1}\right] $ operator $A+\mu $ is positive in $L^{2}\left( G\right) $
for sufficiently large $\mu >0$. Moreover, in view of (1)-(5) all conditons
of Theorem 1 are hold. Then Theorem1 implies the assertion.

\begin{center}
\textbf{9.} \textbf{The Wentzell-Robin type mixed problem for Boussinesq
equations}
\end{center}

Consider the problem $\left( 1.5\right) -\left( 1.6\right) $. \ Let 
\[
\sigma =R^{n}\times \left( 0,1\right) . 
\]

Suppose $\nu =\left( \nu _{1},\nu _{2},...,\nu _{n}\right) $ are nonnegative
real numbers. In this section, we present the following result:

\bigskip \textbf{Theorem 9.1. } Suppose the the following conditions are
satisfied:

(1) $a\left( .\right) \in W_{\infty }^{1}\left( 0,1\right) ,$ $a\left(
.\right) \geq \delta >0,$ $b\left( .\right) ,$ $c\left( .\right) \in
L_{\infty }\left( 0,1\right) ;$

(2) 
\[
\sup\limits_{t\in \left[ 0,1\right] }\left\Vert e^{\left\vert x\right\vert
^{2}\mu ^{2}\left( t\right) }V\left( .,y,t\right) \right\Vert _{L^{\infty
}\left( \sigma \right) }<\infty \text{ for }y\in \left[ 0,1\right] . 
\]

(3) Let $\alpha $, $\beta >0$ and $\alpha \beta <2$. Assume$\ u\in C\left( %
\left[ 0,T\right] ;L^{2}\left( \sigma \right) \right) $ be a solution of the
equation $\left( 1.5\right) -\left( 1.6\right) $ and%
\[
\left\Vert e^{\frac{\left\vert x\right\vert ^{2}}{\beta ^{2}}}u\left(
.,0\right) \right\Vert _{L^{2}\left( \sigma \right) }<\infty ,\left\Vert e^{%
\frac{\left\vert x\right\vert ^{2}}{\alpha ^{2}}}u\left( .,T\right)
\right\Vert _{L^{2}\left( \sigma \right) }<\infty . 
\]

Then $u\left( x,y,t\right) \equiv 0.$

\ \textbf{Proof.} Let us consider operators $A$ in $H=L^{2}\left( G\right) $
that are defined by the equalities 
\[
D\left( A\right) =\left\{ u\in W^{2m,2}\left( G\right) \text{, }B_{j}u=0,%
\text{ }j=1,2,...,m\text{ }\right\} ,\ Au=\sum\limits_{\left\vert \alpha
\right\vert \leq 2m}a_{\alpha }\left( y\right) D_{y}^{\alpha }u\left(
y\right) . 
\]

Then the problem $\left( 8.1\right) -\left( 8.2\right) $ can be rewritten as
the problem $\left( 1.1\right) $, where $u\left( x\right) =u\left(
x,.\right) ,$ $f\left( x\right) =f\left( x,.\right) $,\ $x\in \sigma $ are
the functions with values in\ $H=L^{2}\left( G\right) $. By virtue of $\left[
\text{10, 11}\right] $ the operator $A$ generates analytic semigroup in $%
L^{2}\left( 0,1\right) $. Hence, by virtue of (1)-(5) all conditons of
Theorem 1 are satisfied. Then Theorem1 implies the assertion.

\textbf{References}\ \ 

\begin{enumerate}
\item Agmon S., On the eigenfunctions and on the eigenvalues of general
elliptic boundary value problems, Comm. Pure Appl. Math., 1962, 15,
119-147.\ 

\item H. Amann, Linear and quasi-linear equations,1, Birkhauser, Basel 1995.

\item A. Bonami, B. Demange, A survey on uncertainty principles related to
quadratic forms. Collect. Math. 2006, Vol. Extra, 1--36.

\item C. E. Kenig, G. Ponce, L. Vega, On unique continuation for nonlinear
Schr\"{o}dinger equations, Comm. Pure Appl. Math. 60 (2002) 1247--1262.

\item L. Escauriaza, C. E. Kenig, G. Ponce, L. Vega, On uniqueness
properties of solutions of Schr\"{o}dinger Equations, Comm. PDE. 31, 12
(2006) 1811--1823.

\item L. Escauriaza, C. E. Kenig, G. Ponce, L. Vega, On uniqueness
properties of solutions of the k-generalized KdV, J. of Funct. Anal. 244, 2
(2007) 504--535.

\item L. Escauriaza, C. E. Kenig, G. Ponce, and L. Vega, Hardy's uncertainty
principle, convexity and Schr\"{o}dinger Evolutions, J. European Math. Soc.
10, 4 (2008) 883--907.

\item L. H\"{o}rmander, A uniqueness theorem of Beurling for Fourier
transform pairs, Ark. Mat. 29, 2 (1991) 237--240.

\item J. A. Goldstain, Semigroups of Linear Operators and Applications,
Oxford University Press, Oxfard, 1985.

\item A. Favini, G. R. Goldstein, J. A. Goldstein and S. Romanelli,
Degenerate Second Order Differential Operators Generating Analytic
Semigroups in $L_{p}$ and $W^{1,p}$, Math. Nachr. 238 (2002), 78 --102.

\item V. Keyantuo, M. Warma, The wave equation with Wentzell--Robin boundary
conditions on Lp-spaces, J. Differential Equations 229 (2006) 680--697.

\item Krein S. G., Linear Differential Equations in Banach
space\textquotedblright , American Mathematical Society, Providence, 1971.

\item Lunardi A., Analytic Semigroups and Optimal Regularity in Parabolic
Problems, Birkhauser, 2003.

\item Lions, J-L., Magenes, E., Nonhomogenous Boundary Value Broblems, Mir,
Moscow, 1971.

\item V. B. Shakhmurov, Nonlinear abstract boundary value problems in
vector-valued function spaces and applications, Nonlinear Anal-Theor., v.
67(3) 2006, 745-762.

\item V. B. Shakhmurov, Coercive boundary value problems for regular
degenerate differential-operator equations, J. Math. Anal. Appl., 292 ( 2),
(2004), 605-620.

\item P. Guidotti, \ Optimal regularity for a class of singular abstract
parabolic equations, J. Differential Equations, v. 232, 2007, 468--486.

\item R. Shahmurov, On strong solutions of a Robin problem modeling heat
conduction in materials with corroded boundary, Nonlinear Anal., Real World
Appl., v.13, (1), 2011, 441-451.

\item R. Shahmurov, Solution of the Dirichlet and Neumann problems for a
modified Helmholtz equation in Besov spaces on an annuals, J. Differential
equations, v. 249(3), 2010, 526-550.

\item \ E. M. Stein, R. Shakarchi, Princeton, Lecture in Analysis II.
Complex Analysis, Princeton, University Press (2003).

\item A. Sitaram, M. Sundari, S. Thangavelu, Uncertainty principles on
certain Lie groups, Proc. Indian Acad. Sci. Math. Sci. 105 (1995), 135-151.

\item L. Weis, Operator-valued Fourier multiplier theorems and maximal $%
L_{p} $ regularity, Math. Ann. 319, (2001), 735-758.

\item H. Triebel, Interpolation theory, Function spaces, Differential
operators, North-Holland, Amsterdam, 1978.

\item S. Yakubov and Ya. Yakubov, \textquotedblright Differential-operator
Equations. Ordinary and Partial \ Differential Equations \textquotedblright
, Chapman and Hall /CRC, Boca Raton, 2000.

\item Y. Xiao and Z. Xin, On the vanishing viscosity limit for the 3D
Navier-Stokes equations with a slip boundary condition. Comm. Pure Appl.
Math. 60, 7 (2007), 1027--1055.
\end{enumerate}

\end{document}